\documentclass[a4paper, 12pt]{amsart}

\usepackage{amssymb,amsmath,amsthm,amscd}

\usepackage{amscd}
\usepackage{mathrsfs}
\usepackage{color}
\usepackage[all]{xy}

\def\today{\number\day\space\ifcase\month\or   January\or February\or
   March\or April\or May\or June\or   July\or August\or September\or
   October\or November\or December\fi\   \number\year}

\theoremstyle{definition}
\newtheorem{thm}{Theorem}[section]
\newtheorem{lem}[thm]{Lemma}
\newtheorem{prp}[thm]{Proposition}
\newtheorem{dfn}[thm]{Definition}
\newtheorem{cor}[thm]{Corollary}

\newtheorem{rmk}[thm]{Remark}

\newtheorem{exa}[thm]{Example}

\newtheorem{qst}{Question}

\newcommand{\beq}{\begin{equation}}

\newcommand{\eeq}{\end{equation}}

\newcommand{\beqr}{\begin{eqnarray*}}

\newcommand{\eeqr}{\end{eqnarray*}}

\newcommand{\bal}{\begin{align*}}

\newcommand{\eal}{\end{align*}}

\newcommand{\bei}{\begin{itemize}}

\newcommand{\eei}{\end{itemize}}

\newcommand{\Z}{{\mathbb{Z}}}

\newcommand{\R}{{\mathbb{R}}}

\newcommand{\C}{{\mathbb{C}}}

\newcommand{\N}{{\mathbb{N}}}

\newcommand{\Tr}{{\mathrm{Tr}}}

\title[Limit of iteration of the induced Aluthge transformations]
{Limit of iteration of the induced Aluthge transformations
of centered operators}

\author[H. Osaka]{Hiroyuki Osaka}
\address{ Department of Mathematical Sciences\\ Ritsumeikan University\\ Kusatsu, Shiga, 525-8577  Japan}
\email[]{osaka@se.ritsumei.ac.jp}

\author[T. Yamazaki]{Takeaki Yamazaki}
\address{Department of Electrical, Electronic and Communications Engineering, Toyo University, Kawagoe-Shi, Saitama, 350-8585, Japan}
\email[]{t-yamazaki@toyo.jp}

\date{\today}

\thanks{$^*$Research of the first author is partially supported by the JSPS grant for Scientific Research No.20K03644}



\begin{document}

\begin{abstract}
Aluthge transform is a well-known mapping defined on bounded linear operators. Especially, the convergence property of its iteration has been studied by many authors. In this paper, we discuss the problem for the induced Aluthge transforms which is a generalization of the Aluthge transform defined in 2021. We  give the polar decomposition of the induced Aluthge transformations of centered operators and show its iteration 
converges to a normal operator. In particular, if $T$ is an invertible centered matrix, then iteration of any induced Aluthge transformations converges. Using the canonical standard form of matrix algebras we show that the iteration of any induced Aluthge transformations with respect to the weighted arithmetic mean and the power mean converge. Those observation are extended to the $C^*$-algebra of compact operators on an infinite dimensional Hilbert space, and as an application we show the stability of $\mathcal{AN}$ and $\mathcal{AM}$ properties under the iteration of the induced Aluthge transformations. We also provide concrete forms of their limit points
for centered matrices and several examples. Moreover, we discuss the limit point of the induced Aluthge transformation with respect to the power mean in the injective $II_1$-factor $\mathcal{M}$ and determine the form of its limit for some centered operators in $\mathcal{M}$.
\end{abstract}

\maketitle

\tableofcontents

%\newpage
%%%%%%%%%%%%%%%%%%%%%%%%%%%%%%%%%%%%%%%%%%%%%%%%%%%%%%%%%%%%%%%%%%%%%%%
\section{Introduction}

Let $H$ be a complex Hilbert space, and 
$\mathcal{B}(H)$ be a $C^{*}$-algebra of
bounded linear operators on $H$.
$\mathcal{S}_{2}(H)\subset \mathcal{B}(H)$ 
denotes
the set of all Hilbert-Schmidt operators, 
and let $\|X\|_{2}:=(Tr(X^*X))^{1/2}$ be the
Schatten $2$-norm for $X\in \mathcal{S}_{2}(H)$.
Let $\mathcal{M}_{m}(\mathbb{C})$ be 
the set of all $m$--by--$m$ matrices with 
complex entries.
Let $\mathcal{F}(H)$ and $\mathcal{K}(H)$
be the set of all finite rank operators and
the set of all compact operators on 
$H$, respectively.
For a $C^{*}$-algebra $M$, $M_{+}$ and 
$M_{++}$ are defined as the set of all
positive elements in $M$, and 
the set of all positive invertible 
elements in $M$, respectively.
For self-adjoint $A,B\in \mathcal{B}(H)$, $A\leq B$ (resp. $A<B$) is defined by $B-A\in \mathcal{B}(H)_{+}$ (resp. $B-A\in \mathcal{B}(H)_{++}$). In what follows, the ``$\lim$'' means the strong operator topology. A continuous function $f$ defined on an interval $I\subset \mathbb{R}$ is called operator monotone if $f(A)\leq f(B)$ holds for all self-adjoint $A,B\in \mathcal{B}(H)$ satisfying $A \leq B$ and $\sigma(A),\sigma(B)\subseteq I$, where $\sigma(X)$ means the spectrum of $X\in \mathcal{B}(H)$.

Let $T = U|T| \in \mathcal{B}(H)$ be the polar decomposition of $T$. The {\it Aluthge transformation} \cite{Aluthge 1990} $\Delta(T)$ is defined as follows:
$$
\Delta(T) = |T|^\frac{1}{2}U|T|^\frac{1}{2}.
$$
The Aluthge transform has a lot of nice properties. Especially, the following properties are well known: (i) $\sigma(\Delta(T)) = \sigma(T)$ \cite{Huruya 1997}, (ii) $\Delta(T)$ has a non-trivial subspace if and only if so does $T$ \cite{JKP 2000}, (iii) if $T$ is semi-hyponormal (i.e. $|T^*| \leq |T|$), then $\Delta(T)$ is hyponormal (i.e., $|\Delta(T)^*|^2\leq |\Delta(T)|^2$) \cite{Aluthge 1990}. More generally, for any scalar $\lambda \in [0, 1]$, the $\lambda$-Aluthge transformation of $T$ is defined by $\Delta_\lambda(T):= |T|^{1-\lambda} U|T|^\lambda$ in \cite{Huruya 1997}.
There are many papers on Aluthge transform, for example, \cite{Aluthge 1990, A2004,  AY 2003, APS2011, CJL2005, DS 2009, Furutabook, GOUY pre2, HT 2007, HT 2010, Huruya 1997, JKP 2000, Y 2002}.
In \cite{LLY2014}, as another operator mapping, {\it mean transform} $\hat{T}$ is defined as follows:
$$ \hat{T}:=\frac{U|T|+|T|U}{2}. $$
The mean transform has similar properties to the Aluthge transform, for example, \cite{BCKL2023, CCM2019, CM2020, JKL2020}.
However, $\sigma(\hat{T})=\sigma(T)$ does not hold, in general \cite{LLY2014}.

For $A,B\in \mathcal{B}(H)_{++}$, an operator connection 
$$ A^{\frac{1}{2}}(A^{-\frac{1}{2}}BA^{-\frac{1}{2}})^{\frac{1}{2}}A^{\frac{1}{2}}$$
is called the operator geometric mean. It was first considered by Pusz-Woronowicz in \cite{PW1975}. Then Kubo and Ando generalized it with an axiomatic definition in \cite{KA1980}. Moreover,
for any operator mean $\frak{m}$, there exists a positive 
operator monotone function $f$ defined on $(0,\infty)$ such that 
$f(1)=1$ and 
\begin{equation}
\frak{m}(A,B)=A^{\frac{1}{2}}f(A^{-\frac{1}{2}}BA^{-\frac{1}{2}})
A^{\frac{1}{2}}
\label{eq:operator mean}
\end{equation}
holds for all $A\in \mathcal{B}(H)_{++}$ and $B\in \mathcal{B}(H)_{+}$ in \cite{KA1980}.
In this case, $f$ is called the {\it representing function} of an 
operator mean $\frak{m}$. In what follows, we write $\frak{m}_f$ by 
an operator mean with the representing function $f$.

Recently, one of 
the author defined the induced  Aluthge transform in the viewpoint %of the axiom 
of operator means, which is defined by using double operator integrals in 
\cite{Y 2021}. It interpolates between mean 
and Aluthge transformations when
$|T|$ is invertible.

\begin{dfn}[Induced Aluthge transformation, \cite{Y 2021}]
Let $T = U|T| \in \mathcal{B}(H)$ with the spectral decomposition $|T| = \int_{\sigma(T)}sdE_s$. For an operator mean $\frak{m}_f$ with a 
representing function $f$, the  
{\it induced Aluthge transformation} (IAT, for short)
$\Delta_{\frak{m}_f}(T)$ of  $T$ with respect to $\frak{m}_f$ is defined as follows.
\begin{itemize}
\item[{\rm (i)}]
If $|T|$ is invertible, then 
$$
\Delta_{\frak{m}_f}(T) := \int_{\sigma(|T|)}\int_{\sigma(|T|)}P_f(s, t)dE_s UdE_t,
$$
where $P_f(s, t) = sf(\frac{t}{s})$ for $s, t \in (0, \infty)$.
\item[{\rm (ii)}]
If $|T|$ is not invertible, and if there exists an isometry $V$ such that $T_\varepsilon = V(|T| + \varepsilon I_H)$ is the polar decomposition for all $\varepsilon > 0$ and $\lim_{\varepsilon \searrow 0}T_\varepsilon = T$, then 
$$
\Delta_{\frak{m}_f}(T) = \lim_{\varepsilon\searrow 0}\Delta_{\frak{m}_f}(T_\varepsilon).
$$
\end{itemize}
\end{dfn}

We remark that $\Delta_{\frak{m}_f}(T)$ can be defined 
if $f$ is a representing function of an operator mean \cite{Y 2021}.

We often use the symbol $\Delta_{f}$ instead of $\Delta_{\frak{m}_{f}}$, simply.

\vskip 2mm

\begin{exa}[{\cite[Example 2]{Y 2021}}]\label{ex:induced Aluthge transform}
Let $T \in \mathcal{B}(H)$ such that $|T|$ is 
invertible with
the polar decomposition $T = U|T|$.

\vskip 2mm

\begin{itemize}
\item[{\rm (i)}]
Let $\lambda \in [0, 1]$ and $f_\lambda (x) = 1 - \lambda + \lambda x$ for $x \in [0, \infty)$, i.e., the corresponding operator mean $\frak{m}_{f_\lambda}$ is called the 
$\lambda$-weighted arithmetic mean. 
Then the IAT with respect to $\frak{m}_{f_\lambda}$ is
$$
\Delta_{f_\lambda}(T) = (1 - \lambda)|T|U + \lambda U|T|.
$$
Especially,  $\Delta_{f_{1/2}}(T)=\hat{T}$, the mean transformation
of $T$ defined in \cite{LLY2014}.

\item[{\rm (ii)}]
Let $\lambda \in [0, 1]$ and $g_\lambda (x) = x^\lambda$ for $x \in [0, \infty)$,
i.e.,  the corresponding operator mean $\frak{m}_{g_\lambda}$ is called the 
$\lambda$-weighted geometric mean. Then the IAT with respect to $\frak{m}_{g_\lambda}$ is
$$
\Delta_{g_\lambda} = |T|^{1-\lambda}U|T|^\lambda.
$$
We obtain $\Delta_{g_{1/2}} = \Delta$, 
the Aluthge transform \cite{Aluthge 1990}.

\item[{\rm (iii)}] Let $\lambda\in [0,1]$ and 
$f_{r, \lambda}(x)=[1- \lambda + \lambda x^r]^{\frac{1}{r}}$.
If $r\in [-1,1]\backslash\{0\}$, then $f_{r,\lambda}$ is a 
representing function of the operator power mean.
For a natural number $n$, 
the IAT %induced Aluthge transformation 
with respect to the power mean with 
the parameter $r=\frac{1}{n}$ 
is given as follows.
\begin{equation}
\Delta_{f_{\frac{1}{n}, \lambda}}(T)
 =
\left(
\sum_{k=0}^{n} \begin{pmatrix} n \\ k \end{pmatrix} (1-\lambda)^{k}\lambda^{n-k}|T|^{\frac{k}{n}}U
|T|^{\frac{n-k}{n}} \right).
\label{eq:power mean}
\end{equation}
\end{itemize}
\end{exa}

In this paper, we shall consider iteration of IATs.
For a non-negative integer $n$ and a representing function $f$ of an operator mean, define $\Delta_{f}^{n}(T)$ as follows.
$$ \Delta_{f}^{0}(T)=T \text{ and } 
\Delta_{f}^{n}(T)=\Delta_{f}(\Delta_{f}^{n-1}(T)) $$
for invertible $T\in \mathcal{B}(H)$. We call $\{\Delta^{n}_{f}(T)\}_{n=0}^{\infty}$ the induced Aluthge sequence (the IAS, for short). It is known that 
if $f(x)=\sqrt{x}$ or $f(x)=\frac{1+x}{2}$ and $\dim H<+\infty$, 
then the IAS converges to a normal matrix 
in \cite{APS2011, Y 2021}, respectively.
However concrete forms of their limit points are not known.
Moreover the proof for the case $f(x)=\sqrt{x}$ is complicated.

In this paper, we shall discuss convergence properties of IASs %stability of the iteration of the induced Aluthge transformations 
under the centered condition.
Moreover, we give a concrete limit point. 
As an application
we shall show that 
limit point of IASs stable under  $\mathcal{AN}$ and $\mathcal{AM}$ properties (the definitions of $\mathcal{AN}$ and $\mathcal{AM}$ properties will be introduced in the Section 7). 
Recall that $T\in \mathcal{B}(H)$ is called {\it centered}
if and only if 
$$ \mathcal{O}(T)=\{ U^{m*}|T|U^{m}, |T|, U^{n}|T|U^{n*}; 
n,m=1,2,...\} $$
is a commuting set, i.e., for each two elements in $\mathcal{O}(T)$ are commuting, 
%where $T = U|T| \in \mathcal{B}(H)$ be the polar decomposition 
\cite{MM1974}. 
The set of centered operators includes scalar weighted shifts and isometries. Moreover it is shown that the partial isometry part of the polar decomposition of centered operators is a power partial isometry \cite{MM1974}.

This paper is organized as follows:
In Section 2, we shall give the polar decomposition of IATs
of centered operators. In Section 3, we shall show that the IAS with respect to arbitrary operator mean of an invertible semi-hyponormal centered operator converges to a normal operator. In Section 4, we shall show that if $T$ is an invertible centered matrix, then the IAS with respect to an arbitrary operator mean
converges to a normal matrix. 
Using the canonical standard form of matrix algebras, we show that 
the IAS with respect to the weighted arithmetic mean converges, which is a generalization of \cite[Theorem 5.1]{Y 2021}, and the power mean converges. Moreover it is shown that the limit point does not depend on the weight parameter in Section 5. Those observations can be extended in the $C^*$-algebra of compact operators on an infinite dimensional Hilbert space in Section 6. As applications of the results in Section 6,  we show the stability of $\mathcal{AN}$ and $\mathcal{AM}$ properties under the limit of IAS %iteration of the induced Althuge transformations 
 in Section 7.  
In section 8 we shall discuss the concrete form of the limit point of IAS %$\Delta^n(T)$ 
for a centered matrix and give an example of the limit point. Lastly, we discuss the limit point of the IAS %induced Althuge transformation 
with respect to the power mean in the injective $II_1$-factor and determine the form of its limit which is an generalization of \cite[Theorem 5.2]{DS 2009}. 

%%%%%%%%%%%%%%%%%%%%%%%%%%%%%%%%%%%%%%%%%%%%%%%%%%%%%%%%%%%%%%%%%%%
\section{The induced Aluthge transformation of centered operators}

In this section, we shall give the polar 
decomposition of the IAT %induced Aluthge transformations 
of centered operators.
The following result gives another formula of the IAT.%induced 
%Aluthge transformation.}

%Firstly, we recall the definition of the 
%induced Aluthge transformation. 
%Then we shall give the polar decomposition
%of the induced Aluthge transformations 
%of centered operators.

%\begin{dfn}[The induced Aluthge transformations, {\cite[Definition~4.4]{Y 2021}}]
%\label{defi:generalized Aluthge transformation}
%Let $T=U|T| \in \mathcal{B}(H)$ 
%with the 
%spectral decomposition $|T|=\int_{\sigma(|T|)}sdE_{s}$.
%For an operator mean $\mathfrak{M}_{f}$,
%the  {\it induced Aluthge transformation} 
%$\Delta_{\mathfrak{M}_{f}}(T)$ of $T$
%with respect to $\mathfrak{M}_{f}$
%is defined as follows.
%
%\begin{itemize}
%\item[(1)] If $|T|$ is invertible, then 
%
%$$ \Delta_{\mathfrak{M}_{f}}(T):=
%\int_{\sigma(|T|)}\int_{\sigma(|T|)}\mathcal{P}_{f}(s,t)
%dE_{s}UdE_{t}. $$
%
%\item[(2)] If $|T|$ is not invertible, and
%if there exists an isometry $V$ such that
%$T_{\varepsilon}:=V(|T|+\varepsilon I)$ is 
%the polar decomposition for all $\varepsilon >0$ and 
%$\displaystyle s-\lim_{\varepsilon\searrow 0}T_{\varepsilon}=T$, then
%
%$$ \Delta_{\mathfrak{M}_{f}}(T):= 
%s-\lim_{\varepsilon\searrow 0}
%\Delta_{\mathfrak{M}_{f}}(T_{\varepsilon}). $$
%
%\end{itemize}
%\end{dfn}

\begin{thm}(\cite[Theorem~4.1]{Y 2021})
\label{thm:generalized Aluthge transformation another form}
Let $T \in \mathcal{B}(H)$ 
be invertible and the polar decomposition of $T$ is $T=U|T|$.
For an operator mean $\mathfrak{m}_{f}$,
the  {\it induced Aluthge transformation} 
$\Delta_{f}(T)$ of $T$
with respect to $\mathfrak{m}_{f}$
is computed as follows.
There exists a positive 
probability measure $d\mu$ on $[0,1]$
such that 
$$ \Delta_{f}(T):=
\int_{0}^{1}\int_{0}^{\infty}
e^{-(1-\lambda)s|T|^{-1}}U
e^{-\lambda s|T|^{-1}}dsd\mu(\lambda).
$$
\end{thm}

To consider IAS, %iteration of the induced Aluthge transformations, 
we need the polar decomposition of the IAT.%induced Aluthge transformations. 
The polar decomposition of IATs are obtained in case of the arithmetic and geometric means in \cite{CCM2019} and \cite{IYY2004}, respectively.
Unfortunately, it has not 
been obtained in the cases of arbitrary operator means.
%, generally. 
However we can get 
the polar decomposition of IAT if $T$ is 
an invertible centered operator. 

\begin{dfn}[Centered operators \cite{MM1974}]
\label{dfn:centered}
Let $T=U|T|\in \mathcal{B}(H)$ be 
the polar decomposition. 
\begin{itemize}
\item[(i)] $T$ is called {\it binormal}
if and only if $[ |T|, |T^{*}|]=0$ (i.e., 
$|T| |T^{*}|=|T^{*}| |T|$), and
\item[(ii)] $T$ is called {\it centered}
if and only if 
$$ \mathcal{O}(T)=\{ U^{m*}|T|U^{m}, |T|, U^{n}|T|U^{n*}; 
n,m=1,2,...\} $$
is a commuting set, i.e., 
each pair of elements in $\mathcal{O}(T)$ commutes.
\end{itemize}
\end{dfn}

%The set of all centered operators includes
%all scalar weighted shifts, all normal operators, and all power partial isometries.
%Moreover the structure of centered operators
%is shown in \cite{MM1974}.

\begin{thm}\label{thm:induced Aluthge transform for binormal operators}
Let $T=U|T|\in \mathcal{B}(H)$ be 
the polar decomposition, and $T$ is invertible.
Assume that one of the following conditions are satisfied,
\begin{itemize}
\item[(i)] $\frak{m}$ is a weighted operator arithmetic mean,
\item[(ii)] $\frak{m}$ is an operator mean and $T$ is binormal.
\end{itemize}
Then the polar decomposition of IAT of $T$ with respect to $\frak{m}$ is given as follows.
$$ \Delta_{\frak{m}}(T)=U \frak{m}(U^{*}|T|U, |T|). $$
%
%and it is the polar decomposition.
Hence $|\Delta_{\frak{m}}(T)|=
\frak{m}(U^{*}|T|U,|T|)$.
\end{thm}

To prove this, we refer to the following result.

\begin{thm}[\cite{H 2013}]\label{thm:Hansen}
Let $f$ be a positive operator monotone function defined on the positive half-line
satisfying $f(1)=1$. Then there exists a bounded positive probability measure 
$\mu$ on the closed interval $[0,1]$ 
such that
$$ f(x)=\int_{0}^{1} [1-\lambda+\lambda x^{-1}]^{-1}d\mu(\lambda) \quad
\text{for all $x>0$.} $$
\end{thm}

For $A,B\in \mathcal{B}(H)_{++}$ and any operator mean $\frak{m}$, 
%can be 
%represented by using a positive 
%operator monotone function $f$ such that
%$f(1)=1$ as follows.
%
%$$ \frak{m}(A,B)=A^{1/2}f(A^{-1/2}BA^{-1/2})A^{1/2}. $$
%
%Then using Theorem \ref{thm:Hansen},
we have
\begin{equation}
\frak{m}(A,B)=\int_{0}^{1}[(1-\lambda)A^{-1}+
\lambda B^{-1}]^{-1}d\mu(\lambda)
\label{eq:operator mean_Hansen}
\end{equation}
by \eqref{eq:operator mean} and Theorem \ref{thm:Hansen}.

\begin{proof}[Proof of Theorem \ref{thm:induced Aluthge transform for binormal operators}]
(i) Assume that $\frak{m}$ is a weighted operator arithmetic mean. Then
\begin{align*}
\Delta_{\frak{m}}(T) & = (1-\lambda)|T|U+\lambda U|T| \\
& = 
U [(1-\lambda)U^{*}|T|U+\lambda |T|] \\
& = 
U\frak{m}(U^{*}|T|U,|T|), 
\end{align*}
since $U$ is unitary. $|\Delta_{\frak{m}}(T)|=[(1-\lambda)U^*|T|U+\lambda |T|]$ can be shown easily.

(ii) If $T$ is binormal, then  
\begin{equation}
\begin{split}
[U^{*}|T|U,|T|] & = 
U^{*}|T|U|T|-|T|U^{*}|T|U \\
& =
U^{*}(|T|U|T|U^{*}-U|T|U^{*}|T|)U\\
& =
U^{*}(|T||T^{*}|-|T^{*}||T|)U=0.
\end{split}
\label{eq:binormal}
\end{equation}
We note that $e^{A+B}=e^{A}e^{B}$
holds for self-adjoint $A,B\in \mathcal{B}(H)$
satisfying $[A,B]=0$.
Then we have
\begin{align*}
\Delta_{\frak{m}}(T) 
& =
\int_{0}^{1}\int_{0}^{\infty} 
e^{-(1-\lambda)s|T|^{-1}}U
e^{-\lambda s|T|^{-1}}ds d\mu(\lambda) 
\quad \text{by Theorem \ref{thm:generalized Aluthge transformation another form}}
\\
& =
U \int_{0}^{1}\int_{0}^{\infty} 
U^{*} e^{-(1-\lambda)s|T|^{-1}}U
e^{-\lambda s|T|^{-1}}ds d\mu(\lambda)\\
& \qquad \text{(since $T$ is invertible and $U$ is unitary)}\\
& =
U \int_{0}^{1}\int_{0}^{\infty} 
e^{-(1-\lambda)s(U^{*} |T|U)^{-1}}
e^{-\lambda s|T|^{-1}}ds d\mu(\lambda) 
\\
& =
U \int_{0}^{1}\int_{0}^{\infty} 
e^{-(1-\lambda)s(U^{*} |T|U)^{-1}
-\lambda s|T|^{-1}}ds d\mu(\lambda)
\quad\text{(by \eqref{eq:binormal})}
\\
& =
U \int_{0}^{1} 
\left[(1-\lambda)(U^{*} |T|U)^{-1}
+\lambda |T|^{-1}\right]^{-1} 
d\mu(\lambda)
\\
& =
U \frak{m}(U^{*}|T|U,|T|)
\quad \text{(by \eqref{eq:operator mean_Hansen}).}
\end{align*}
\end{proof}

Next, we shall consider IASs. 
We will show that (i) if $T\in \mathcal{B}(H)$ is invertible and centered, 
then $\Delta_{\frak{m}}(T)$ is so. 
Hence $\Delta^{n}_{\frak{m}}(T)$ is 
centered for all non-negative integer $n$,
and (ii) we will give a concrete formula of
the polar decomposition of 
$\Delta^{n}_{\frak{m}}(T)$ for an invertible centered 
$T\in \mathcal{B}(H)$. To show them, we prepare 
the following lemmas.

\begin{lem}\label{lem commute operator mean}
Let $\frak{m}$ be an operator mean with the
representing function $f$, and
let $\{A,B,C,D\}\subset \mathcal{B}(H)_{++}$ 
be a commuting set. Then 
$$ [\frak{m}(A,B), \frak{m}(C,D)]=0. $$
\end{lem}

\begin{proof}
We note that, if $[X,Y]=0$, then 
$[f(X),f(Y)]=0$. 
Since $\frak{m}(A,B)=A^{1/2}f(A^{-1/2}BA^{-1/2})A^{1/2}$ and 
$\{A,B,C,D\}$ is a commuting set, 
we have $ [\frak{m}(A,B), \frak{m}(C,D)]=0. $
\end{proof}

\begin{lem}\label{lem: centerd}
Let $T\in \mathcal{B}(H)$ be invertible and centered.
Then $\Delta_{\frak{m}}(T)$ is 
invertible and centered.
\end{lem}

\begin{proof}
Since $T$ is centered, $T$ is binormal, and we have the polar decomposition of 
$\Delta_{\frak{m}}(T)$ as
$\Delta_{\frak{m}}(T)=U \frak{m}(U^{*}|T|U, |T|) $
by Theorem \ref{thm:induced Aluthge transform for binormal operators}.
Then for any natural numbers $n$ and $m$, we have
\begin{align*}
& \{|\Delta_{\frak{m}}(T)|, 
U^{*m}|\Delta_{\frak{m}}(T)|U^{m},
U^{n}|\Delta_{\frak{m}}(T)|U^{*n}\}\\
& =
\{\frak{m}(U^{*}|T|U, |T|), 
U^{*m}\frak{m}(U^{*}|T|U, |T|)U^{m},
U^{n}\frak{m}(U^{*}|T|U, |T|)U^{*n}\}\\
& =
\{\frak{m}(U^{*}|T|U, |T|), 
\frak{m}(U^{*(m+1)}|T|U^{m+1}, U^{*m}|T|U^{m}),\\
& \hspace*{5cm}
\frak{m}(U^{n-1}|T|U^{*(n-1)}, U^{n}|T|U^{*n})\}.
\end{align*}
Since $T$ is centered and Lemma \ref{lem commute operator mean}, 
it is a commuting set for all $n,m\in \mathbb{N}$.
Therefore $\Delta_{\frak{m}}(T)$ is centered.
Invertibility of $\Delta_{\frak{m}}(T)$
is obtained by Theorem \ref{thm:induced Aluthge transform for binormal operators}
since $|\Delta_{\frak{m}}(T)|=\frak{m}(U^{*}|T|U,|T|)$ is invertible.
\end{proof}

\begin{thm}\label{thm:iteration}
Let $T=U|T|\in \mathcal{B}(H)$ be 
the polar decomposition, and $T$ is invertible.
Then the following statements hold.
\begin{itemize}
\item[(i)] If $\frak{m}$ is an weighted operator arithmetic mean, 
then for each 
non-negative integer $n$, the polar decomposition of 
$ \Delta_{\frak{m}}^{n}(T)$
is 
\begin{equation}
\Delta_{\frak{m}}^{n}(T)=U \frak{M}_{n}(U, |T|), 
\label{eq:Theorem iteration}
\end{equation}
where $\frak{M}_{0}(U,|T|)=|T|$ and
$$\frak{M}_{n}(U,|T|):=\frak{m}\biggl(U^{*}\frak{M}_{n-1}(U,|T|)U, 
\frak{M}_{n-1}(U,|T|)\biggr)> 0$$ 
for $n=1,2,...$.
\item[(ii)] If $\frak{m}$ is an operator mean 
and $T$ is centered. 
Then for each 
non-negative integer $n$, $ \Delta_{\frak{m}}^{n}(T)$ is centered and
the polar decomposition of 
$ \Delta_{\frak{m}}^{n}(T)$
is given as the form \eqref{eq:Theorem iteration}.
\end{itemize}
\end{thm}

\begin{proof}
(i) can be shown in the same way as (ii). Hence we only prove (ii).
(ii) Firstly, we can obtain that 
$ \Delta_{\frak{m}}^{n}(T)$
is centered by Lemma \ref{lem: centerd}.
Hence we only show \eqref{eq:Theorem iteration}.
It can be shown by mathematical induction 
on $n$. The $n=0$ case is obvious. In fact, 
since $T$ is centered, 
$ \Delta_{\frak{m}}^{0}(T)=T=U|T|$ is centered.
%Hence it is binormal.
Assume that Theorem \ref{thm:iteration} 
holds for all non-negative integer $n$ such that
$0\leq n\leq k$.

We shall prove $n=k+1$ case.
By the induction hypothesis, 
$$ \Delta_{\frak{m}}^{k}(T)=U \frak{M}_{k}(U, |T|) $$
is the polar decomposition, and 
$\Delta^{k}_{\frak{m}}(T)$ is binormal.
Then by Theorem \ref{thm:induced Aluthge transform for binormal operators}, we have
\begin{align*}
\Delta_{\frak{m}}^{k+1}(T) & =
U
\frak{m}(U^{*}|\Delta_{\frak{m}}^{k}(T)|U, |\Delta_{\frak{m}}^{k}(T)|)\\
& =
U
\frak{m}(U^{*}\frak{M}_{k}(U, |T|)U, \frak{M}_{k}(U, |T|))\\
& = U \frak{M}_{k+1}(U,|T|).
\end{align*}
\end{proof}

We can obtain a concrete form of 
$\frak{M}_{k}(U, |T|)$ if $T$ is an 
invertible centered operator and 
$\frak{m}$ is the power mean case.
For $r\in [-1,1]\setminus\{0\}$ and
$\lambda\in [0,1]$, let
$f_{r,\lambda}(x)=[1-\lambda +\lambda x^r]^{1/r}$. Then 
$f_{r,\lambda}$ is an operator monotone function, and the
corresponding operator mean is called 
the power mean. We remark that 
$f_{r,\lambda}(x)$ is monotone increasing on 
$r\in [-1,1]\setminus\{0\}$ and 
$\lim_{r\to 0}f_{r,\lambda}(x)=x^{\lambda}$.

\begin{prp}\label{prp: power mean_centered}
Let $T$ be an invertible centered operator with the polar decomposition $T=U|T|$. For $\lambda\in (0,1)$ and $r\in [-1,1]\setminus\{0\}$. 
%For a non-negative integer $k\in [0,n]$, let 
%$a_{k}:=(1-\lambda)^{k}\lambda^{n-k}$ and 
Let
$f_{r,\lambda}(x)=[1-\lambda+\lambda x^{r}]^{1/r}$ and $\frak{m}$ be an operator mean with a representing function $f_{r,\lambda}$. 
Then
\begin{equation}
\frak{M}_{n}(U,|T|)= \left[ \sum_{k=0}^{n}\begin{pmatrix}n \\ k \end{pmatrix} (1-\lambda)^{k}\lambda^{n-k}U^{*k}|T|^{r}U^{k}\right]^{1/r} 
\label{eq:iteration-power mean}
\end{equation}
for all non-negative integer $n$.
\end{prp}

\begin{rmk}\label{rm:r=1 case}
\begin{itemize}
\item[(i)] If $r=1$, then the centered condition of $T$ can be removed by Theorem 
\ref{thm:iteration} (i).
\item[(ii)] Proposition \ref{prp: power mean_centered} is an extension of \cite{CCM2019}
which showed the case $r=1$ and $\lambda=\frac{1}{2}$ without the centered condition.
\end{itemize}
\end{rmk}

\begin{proof}[Proof of Proposition~\ref{prp: power mean_centered}]

We note that for $A,B\in \mathcal{B}(H)_{++}$ such that $[A,B]=0$, the following holds.
\begin{align*}
\frak{m}(A,B) & =
A^{1/2}\left[(1-\lambda)I+\lambda (A^{-1/2}BA^{-1/2})^{r}\right]^{1/r}A^{1/2} \\
& =
\left[(1-\lambda)A^{r}+\lambda B^{r}\right]^{1/r}.
\end{align*}
We shall prove \eqref{eq:iteration-power mean} by mathematical induction on $n$. The case $n=0$. It is easy to see $\frak{M}_{0}(U,|T|)=|T|$. Assume that \eqref{eq:iteration-power mean} holds for some non-negative integers $n$, and we shall prove $n+1$ case. Since $T$ is centered, $U\frak{M}_{n}(U,|T|)$ is centered and it is the polar decomposition for all $n=0,1,2,...$ by Theorem \ref{thm:iteration}.
Moreover
\begin{align*}
\frak{M}_{n+1}(U,|T|)^{r}
& =
\frak{m}\left(U^{*}\frak{M}_{n}(U,|T|)U, \frak{M}_{n}(U,|T|)\right)^{r} \\
& =
(1-\lambda) U^{*}\frak{M}_{n}(U,|T|)^{r}U+\lambda \frak{M}_{n}(U,|T|)^{r} \\
& =
(1-\lambda) U^{*}\left[ \sum_{k=0}^{n}\begin{pmatrix}n \\ k \end{pmatrix} (1-\lambda)^{k}\lambda^{n-k}U^{*k}|T|^{r}U^{k}\right]U \\
& \qquad\qquad +
\lambda \left[ \sum_{k=0}^{n}\begin{pmatrix}n \\ k \end{pmatrix} (1-\lambda)^{k}\lambda^{n-k}U^{*k}|T|^{r}U^{k}\right] \\
& =
\sum_{k=0}^{n}\begin{pmatrix}n \\ k \end{pmatrix} (1-\lambda)^{k+1}\lambda^{n-k}U^{*k+1}|T|^{r}U^{k+1} \\
& \qquad\qquad +
\sum_{k=0}^{n}\begin{pmatrix}n \\ k \end{pmatrix} (1-\lambda)^{k}\lambda^{n+1-k}U^{*k}|T|^{r}U^{k} \\
& =
\sum_{k=0}^{n+1}\begin{pmatrix}n+1 \\ k \end{pmatrix} (1-\lambda)^{k}\lambda^{n+1-k}U^{*k}|T|^{r}U^{k}.
\end{align*}
The proof is completed.
\end{proof}

\begin{rmk}
$f_{0,\lambda}(x)=\lim_{r\to 0}f_{r,\lambda}(x)=x^{\lambda}$ case, we have
\begin{equation}
\frak{M}_{n}(U,|T|)= \prod_{k=0}^{n}
U^{*k}|T|^{\begin{pmatrix}n \\ k \end{pmatrix} 
(1-\lambda)^{k}\lambda^{n-k}}U^{k} 
\label{eq:iteration-geometric mean}
\end{equation}
for all non-negative integer $n$.
\end{rmk}

\medskip

It can be shown by the following fact \cite[Proof of Lemma 3.1]{FFS2009}:
Let $n$ be a non-negative integer, and $A_{k}\in \mathcal{B}(H)_{++}$ $(k=1,2,...,n)$. 
For $a_{k}\in (0,1)$ such that $\sum_{k=1}^{n}a_{k}=1$, then 
$$ \lim_{r\to 0} \left(\sum_{k=1}^{n}a_{k}A_{k}^{r}\right)^{\frac{1}{r}}=
\exp \left(\sum_{k=1}^{n}a_{k}\log A_{k}\right) $$
holds, uniformly.
Put $a_{k}=\begin{pmatrix}n \\ k \end{pmatrix}(1-\lambda)^{k}\lambda^{n-k}$
and $A_{k}=U^{*k}|T|U^{k}$. Then we have \eqref{eq:iteration-geometric mean}
because $T$ is centered.

%%%%%%%%%%%%%%%%%%%%%%%%%%%%%%%%%%%%%%%%%%%%%%%%%%%%%%%%%%%%%%%%%%%%%%%%%%%%%%
\section{Semi-hyponormal centered operators case}

In this section, we shall show that 
IAS of an invertible semi-hyponormal 
centered operator converges to a 
normal operator. An operator $T\in \mathcal{B}(H)$ is said to be semi-hyponormal if and only if $|T^*|\leq |T|$ holds.
In \cite{GOUY pre2}, we showed that 
if $T$ is an invertible semi-hyponormal operator and 
$\frak{m}$ is a non-weighted arithmetic mean,
then the IAS %iteration of the induced Aluthge transformations 
of $T$ converges to a 
normal operator. 
However, we could not have any result when 
$\frak{m}$ is an arbitrary operator mean.
In this section, we discuss the problem 
when $T$ is a centered operator.

Because an arbitrary operator mean can be 
represented via the harmonic mean
by Theorem \ref{thm:Hansen} and 
\eqref{eq:operator mean_Hansen}, 
we first consider the harmonic mean case before proving a general result. 

\begin{prp}\label{prp: harmonic mean}
Let $T=U|T|$ 
be the polar decomposition of an 
invertible centered operator $T\in \mathcal{B}(H)$, $\lambda\in (0,1)$, and $\frak{h}$ be an operator weighted harmonic mean, i.e., the representing function $f$ is $f(x)=[1-\lambda+\lambda x^{-1}]^{-1}$.
If $T$ is semi-hyponormal, then
$\{\Delta^{n}_{\frak{h}}(T)\}$ converges to a normal operator $UL$,
where
$L:=\lim_{n\to \infty}{U^{*}}^{n}|T|U^{n}$. 
\end{prp}

We remark that if $T=U|T|$ is a semi-hyponormal operator, then 
$U|T|U^{*}=|T^{*}|\leq |T|$ holds. Hence we have
$$ |T|\leq U^{*}|T|U\leq U^{*2}|T|U^{2}\leq \cdots \leq 
U^{*n}|T|U^{n}\leq \cdots \leq \|T\|I. $$
Therefore there exists $L=\lim_{n\to \infty}{U^{*}}^{n}|T|U^{n}$.
$L$ is called the polar symbol \cite[page 28]{X1983}.
To prove Proposition \ref{prp: harmonic mean}, we prepare the following lemma.

\begin{lem}\label{lem:limit}
Let $T=U|T|$ be the polar decomposition of an invertible operator $T\in \mathcal{B}(H)$. For $\lambda\in (0,1)$, 
the following hold.
\begin{itemize}
\item[(1)] If $T$ is semi-hyponormal,
then there exists $L=\lim_{n\to \infty}
U^{*n}|T|U^{n}$ and 
$$ \sum_{k=0}^{n} \begin{pmatrix} n \\ k \end{pmatrix} (1-\lambda)^{k}\lambda^{n-k} U^{*k}|T|U^{k} \nearrow L\quad \text{(as $n\to \infty$)}. $$
\item[(2)] If $T^{*}$ is semi-hyponormal,
then there exists $L=\lim_{n\to \infty}
U^{*n}|T|U^{n}$ and 
$$ \sum_{k=0}^{n} \begin{pmatrix} n \\ k \end{pmatrix} (1-\lambda)^{k}\lambda^{n-k}U^{*k}|T|U^{k} \searrow L\quad \text{(as $n\to \infty$)}. $$
\end{itemize}
\end{lem}

\begin{proof}
Proof of (1) is given in the similar way as 
\cite[Theorem 3.4]{GOUY pre2} (also 
see \cite{CCM2019}). 

Proof of (2). If $T^{*}$ is semi-hyponormal 
if and only if $|T|\leq |T^{*}|=U|T|U^{*}$, 
then we have
$$ |T|\geq U^{*}|T|U\geq U^{*2}|T|U^{2}
\geq \cdots \geq U^{*n}|T|U^{n}\geq 0.$$
Hence there exists 
$L=\lim_{n\to \infty} U^{*n}|T|U^{n}$.
The latter part can be proven in the similar way as (1).
\end{proof}

\begin{proof}[Proof of Proposition \ref{prp: harmonic mean}]
Since $T$ is centered and by 
Theorem \ref{thm:induced Aluthge transform for binormal operators}, we have
\begin{align*}
\Delta_{\frak{h}}(T) & =
U \frak{h}(U^{*}|T|U, |T|) \\
& =
U \left[(1-\lambda) (U^{*}|T|U)^{-1}+\lambda |T|^{-1}\right]^{-1}\\
& =
U \left[(1-\lambda) U^{*}|(T^{*})^{-1}|U+\lambda |(T^{*})^{-1}|\right]^{-1}.
\end{align*}
Hence 
\begin{equation}
|\Delta_{\frak{h}}(T)|^{-1}=
(1-\lambda) U^{*}|(T^{*})^{-1}|U+\lambda |(T^{*})^{-1}|.
\label{eq:equality}
\end{equation}
On the other hand, we have the polar decomposition $(T^{*})^{-1}=U|(T^{*})^{-1}|$,
and $(T^{*})^{-1}$ is centered 
since $|(T^{*})^{-1}|=|T|^{-1}$.
Moreover 
let $\frak{a}$ be a weighted 
operator arithmetic mean.
We shall show 
$ |\Delta_{\frak{a}}^{n}((T^{*})^{-1})|=
|\Delta_{\frak{h}}^{n}(T)|^{-1}$ holds for 
all $n=0,1,2,...$ by mathematical induction
on $n$.
The $n=0$ case has already been shown.
We shall show the $n=1$ case.
By \eqref{eq:equality},
\begin{align*}
\Delta_{\frak{a}}((T^{*})^{-1})
& =
U \frak{a}(U^{*}|(T^{*})^{-1}|U, |(T^{*})^{-1}|)\\
& =
U \left[
(1-\lambda) U^{*}|(T^{*})^{-1}|U+\lambda |(T^{*})^{-1}|
\right]
 =
U|\Delta_{\frak{h}}(T)|^{-1},
\end{align*}
we have $|\Delta_{a}((T^{*})^{-1})|=
|\Delta_{h}(T)|^{-1}$.
Assume that it  holds for $n=k$.
Then
\begin{align*}
|\Delta_{\frak{a}}^{k+1}((T^{*})^{-1})|
& =
(1-\lambda) U^{*}|\Delta_{\frak{a}}^{k}((T^{*})^{-1})|U+ \lambda 
|\Delta_{\frak{a}}^{k}((T^{*})^{-1})|\\
& =
(1-\lambda) U^{*}|\Delta_{\frak{h}}^{k}(T))|^{-1}U+ \lambda
|\Delta_{\frak{h}}^{k}(T)|^{-1}
 =
|\Delta_{\frak{h}}^{k+1}(T)|^{-1}.
\end{align*}
Hence $ |\Delta_{\frak{a}}^{n}((T^{*})^{-1})|=
|\Delta_{\frak{h}}^{n}(T)|^{-1}$ holds for 
all $n=0,1,2,...$.
Therefore by Proposition \ref{prp: power mean_centered} for $r=1$,
%
%the similar way to 
%\cite[Theorem 2.15]{CCM2019},
we have
$$
|\Delta_{\frak{a}}^{n}((T^{*})^{-1})|
 = 
\sum_{k=0}^{n} \begin{pmatrix} n \\ k \end{pmatrix} \lambda^{k}(1-\lambda)^{n-k}U^{*k}|((T^{*})^{-1}|U^{k}.
$$

On the other hand,
since $T$ is semi-hyponormal, 
we have
\begin{align*}
(T^{*})^{-1} [(T^{*})^{-1}]^{*} 
& =
(T^{*})^{-1} T^{-1}\\
& =
(TT^{*})^{-1} \\
& 
\geq  
(T^{*}T)^{-1} 
 =
[(T^{*})^{-1}]^{*}(T^{*})^{-1},
\end{align*}
and hence $[ (T^{*})^{-1}]^{*}$ is 
semi-hyponormal. 
Hence by Lemma \ref{lem:limit} (2),
$\{|\Delta_{\frak{a}}^{n}((T^{*})^{-1})|\}$ converges to 
\begin{align*}
\lim_{n\to \infty}U^{*n}|(T^{*})^{-1}|U^{n}
& =
\lim_{n\to \infty}U^{*n}|T|^{-1}U^{n}\\
& =
\lim_{n\to \infty}(U^{*n}|T|U^{n})^{-1}=L^{-1}.
\end{align*}
Therefore 
$$ \lim_{n\to\infty}\Delta_{\frak{h}}^{n}(T)=
U \lim_{n\to\infty}|\Delta_{\frak{h}}^{n}(T)|=
U \lim_{n\to\infty}|\Delta_{\frak{a}}^{n}\left((T^{*})^{-1}\right)|^{-1}=U(L^{-1})^{-1}=UL.
$$
Normality of $UL$ can be shown as follows: Since $U$ is unitary, $UL=LU$ holds, i.e., $UL$ is quasinormal. Moreover since $T$ is invertible, $UL$ is also invertible, and hence $UL$ is normal
(see {\cite[Corollary 1.3]{X1983}}).
\end{proof}

\begin{thm}\label{thm:iteration general}
Let $T=U|T|$ 
be the polar decomposition of an 
invertible centered operator $T\in \mathcal{B}(H)$, and $\frak{m}_f$ be
an operator mean, s.t., 
\begin{equation}
[1-\lambda+\lambda t^{-1}]^{-1}
\leq f(t) \leq 1-\lambda+\lambda t
\label{eq:assumption}
\end{equation}
holds for all $t>0$, where $f$ is a 
representing function 
of $\frak{m}_f$ and $\lambda=f'(1)$.
If $T$ is semi-hyponormal, 
then
$$\lim_{n\to \infty} \Delta^{n}_{f}(T)=
UL. $$
\end{thm}

\begin{proof}
By Proposition \ref{prp: harmonic mean},
$\lim_{n\to \infty} |\Delta_{\frak{a}}^{n}(T)|=
\lim_{n\to \infty} |\Delta_{\frak{h}}^{n}(T)|=L$ holds.
%By \eqref{eq:assumption}, we have
%
%$$ |\Delta_{\frak{h}}^{n}(T)|\leq |%\Delta_{\frak{m}}^{n}(T)|
%\leq |\Delta_{\frak{a}}^{n}(T)|\leq L.$$
%
Then by Theorem \ref{thm:induced Aluthge transform for binormal operators} and \eqref{eq:assumption}, we have
\begin{align*}
|\Delta_{\frak{h}}(T)| 
& =
\frak{h}(U^{*}|T|U,|T|) \\
& \leq 
\frak{m}_f(U^{*}|T|U,|T|)=|\Delta_{f}(T)| \\
& \leq 
\frak{a}(U^{*}|T|U,|T|)=|\Delta_{\frak{a}}(T)|.
\end{align*}
Moreover, by the monotonicity of operator means, 
\begin{align*}
|\Delta^{2}_{\frak{h}}(T)| 
& =
\frak{h}(U^{*}|\Delta_{\frak{h}}(T)|U,|\Delta_{\frak{h}}(T)|) \\
& \leq 
\frak{m}_f(U^{*}|\Delta_{\frak{h}}(T)|U,|\Delta_{\frak{h}}(T)|) \\
& \leq 
\frak{m}_f(U^{*}|\Delta_{f}(T)|U,|\Delta_{f}(T)|) =|\Delta^{2}_{f}(T)| \\
& \leq 
\frak{a}(U^{*}|\Delta_{f}(T)|U,|\Delta_{f}(T)|)\\
& \leq 
\frak{a}(U^{*}|\Delta_{\frak{a}}(T)|U,|\Delta_{\frak{a}}(T)|)=|\Delta^{2}_{\frak{a}}(T)|.
\end{align*}
Repeating the above, we have 
$$ |\Delta_{\frak{h}}^{n}(T)|\leq |\Delta_{f}^{n}(T)|
\leq |\Delta_{\frak{a}}^{n}(T)|\leq L$$
for all $n=1,2,...$, where the last inequality holds by Lemma \ref{lem:limit} (1). 
Moreover $\lim_{n\to \infty}\langle 
|\Delta_{f}^{n}(T)|x,x\rangle =\langle Lx,x\rangle$ for all $x\in H$
by Proposition \ref{prp: harmonic mean}.
Therefore $L-|\Delta_{f}^{n}(T)|\geq 0$ 
and
$$ \|(L-|\Delta_{f}^{n}(T)|)^{1/2}x\|^{2}=
\langle Lx,x\rangle-
\langle |\Delta_{f}^{n}(T)|x,x\rangle\to 0, 
\quad\text{as $n\to \infty$}
$$
for all $x\in H$.
Hence $\lim_{n\to \infty}|\Delta_{f}^{n}(T)|=L$.
Therefore, 
$\lim_{n\to \infty}\Delta^{n}_{f}(T)=UL$.
\end{proof}

\begin{rmk}
\begin{itemize}
\item[(i)] If $f(t)=1-\lambda+\lambda t$ in Theorem \ref{thm:iteration general}.
We do not need the centered condition of $T$ (see \cite[Theorem 3.4]{GOUY pre2}).
(ii) We cannot remove the 
semi-hyponormality of $T$ in Theorem 
\ref{thm:iteration general}.
In fact, there exists a weighted shift 
(weighted shift is centered) such 
that $\{\Delta_{\frak{m}}^{n}(T)\}$ does not converge
in \cite[Theorem 5.2]{Y 2021}.
\end{itemize}
\end{rmk}

\begin{qst}
What is the equivalent condition 
of $T\in \mathcal{B}(H)$
for which 
$\{\Delta_{\frak{m}}^{n}(T)\}$ converges as
$n\to \infty$.
\end{qst}

%If $T\in \mathcal{AN}(H)$ is semi-hyponormal, then $UL\in \mathcal{AN}(H)$ \cite[Theorem 3.5]{GOUY pre2}. 
%Hence 
%$\lim_{n\to \infty}\Delta^{n}_{\frak{m}}(T)\in
%\mathcal{AN}(H)$.
%
%
%\begin{qst}
%Are the following questions hold for 
%centered or non-centered operators?
%
%\begin{itemize}
%\item[(i)] Is there exists 
%$\lim_{n\to \infty}|\Delta_{m}^{n}(T)|$
%for any non semi-hyponormal 
%$\mathcal{AN}$-operator $T$?
%
%\item[(ii)] Is there exists 
%$\lim_{n\to \infty}|\Delta_{m}^{n}(T)|$
%for any compact operator $T$?
%
%\item[(iii)] Is there exists 
%$\lim_{n\to \infty}|\Delta_{m}^{n}(T)|$
%for any $d$--by--$d$ matrix $T$?
%\end{itemize}
%\end{qst}
%
%For the above (ii) and (iii), 
%if $T$ is compact and
%semi-hyponormal, then $T$ is normal.
%Hence $T$ is centered and $T=\Delta_{\frak{m}}(T)$ holds.
%Therefore we only consider the non semi-hyponormal case.
%
%\begin{qst}
%Can we discuss 
%Theorem \ref{thm:iteration general}
%in the setting of \cite{DS 2009}?
%\end{qst}
%
%The second author thinks that 
%\cite[Theorem 5.2]{DS 2009}
%may be extended to the power mean
%case, i.e., the representing function is
%$f(t)=(\frac{1+t^{r}}{2})^{1/r}$.
%And the limit may be 
%$(E^{\alpha}(|T|^{r}))^{1/r}$.

%%%%%%%%%%%%%%%%%%%%%%%%%%%%%%%%%%%%%%%%%%%%%%%%%%%%%%%%%%%%%%%%%%%%%%%
\section{Centered matrices case}

In the previous section, we considered 
semi-hyponormal operators. This condition 
cannot be removed. However, if $T$ is a 
centered matrix, then the IAS
%
%iteration of the induced 
%Aluthge transformations 
converges.
In this section, we shall show the fact.
We note that, if $T$ is a matrix and 
there exists a limit point of 
$\{\Delta_{\frak{m}_{f}}^{n}(T)\}$, then 
the limit point is a normal matrix 
\cite[Theorem 2.4]{Y 2021}.

\begin{thm}\label{thm:matrices case}
Let $T=U|T|\in \mathcal{M}_{m}(\mathbb{C})$ be 
invertible and centered, and let
$\frak{m}_{f}$ be an operator mean
with a representing function $f$ satisfying
$f(x)\leq 1-\lambda +\lambda x$ for all $x>0$, where $\lambda :=f'(1)\in (0,1)$.
Then 
$\{\Delta_{\frak{m}_{f}}^{n}(T)\}$ converges 
to a normal matrix 
as $n\to \infty$.
\end{thm}

\begin{rmk}
The case $f(x)=\frac{1+x}{2}$ has been already shown in 
\cite[Theorem 5.1]{Y 2021} without the centered condition of $T$.
By the similar way, the case $f(x)=1-\lambda+\lambda x$ for $\lambda\in (0,1)$ in 
Theorem \ref{thm:matrices case} also holds without the centered condition of $T$, too.
\end{rmk}

To prove Theorem \ref{thm:matrices case},
we prepare the following results.

\begin{lem}\label{lem:centered matrices}
Let $T=U|T|\in \mathcal{M}_{m}(\mathbb{C})$ be 
the polar decomposition of an 
invertible centered matrix.
Then there exists a non-negative integer
$n$ such that
$|T|=U^{*n}|T|U^{n}$.
\end{lem}

\begin{proof}
Since $U$ is unitary, 
$U^{*k}|T|U^{k}$ has the same 
characteristic function for all 
$k=0,1,2,...$, that is,
every $\sigma(U^{*k}|T|U^{k}) $ 
is the same with the same multiplicity for all $k=0,1,2,...$.
If $T$ is centered, then 
$\{U^{*k}|T|U^{k}: k=0,1,2,...\}$ is commuting.
Hence all $U^{*k}|T|U^{k}$ are 
simultaneous diagonalizable.
Therefore we may assume $|T|=diag(a_{1},...,a_{m})$, and 
every $U^{*k}|T|U^{k}$ is also a diagonal matrix
whose diagonal entries are a permutation 
of $\{a_{1},...,a_{m}\}$. The number of permutations is at most $m!$ This means that there exists $n$ such that $n\leq m!$ and $|T|=U^{*n}|T|U^{n}$. 
\end{proof}

\begin{thm}\label{thm:multivariate operator mean}
Let $A_{1},...,A_{n}\in \mathcal{S}_{2}(H)_{++}$
%be positive invertible operators.
Let $\frak{m}_{f}$ be an 
operator mean with a representing function $f$
such that $f'(1):=\lambda\in (0,1)$ and $f(x)\leq 1-\lambda+\lambda x$
for all $x>0$.
Define $A^{(1)}_{k}=A_{k}$, 
$A_{k}^{(\ell+1)}:=
\frak{m}_{f}(A_{k+1}^{(\ell)}, A_{k}^{(\ell)})$ 
for $k=1,2,...,n-1$ and 
$A_{n}^{(\ell+1)}:=\frak{m}_{f}(A_{1}^{(\ell)},
A_{n}^{(\ell)})$. Then for all $k=1,2,...,n$,
the sequences $\{A_{k}^{(\ell)}\}$ 
converge to the same limit as
$\ell\to \infty$. 
 \end{thm}

Theorem \ref{thm:multivariate operator mean} is firstly considered by P\'{a}lfia \cite{Palfia 2011} for $\lambda=\frac{1}{2}$.
The above result is a generalization of 
P\'{a}lfia's result. Proof is just a modification, 
but we introduce it for the reader's convenience
(see also \cite[Theorem 5.18]{HP 2014 book}).

\begin{proof}
From the assumption
\begin{equation}
\frak{m}_{f}(X,Y)\leq (1-\lambda)X+\lambda Y
\text{\quad for all
$X,Y\in \mathcal{B}(H)_{++}$,}
\label{eq:arithmetic-mean inequality}
\end{equation}
we have
$$ \sum_{k=1}^{n}A_{k}^{(\ell+1)}\leq 
\sum_{k=1}^{n}A_{k}^{(\ell)}.$$
Therefore the decreasing positive sequence 
has a limit $L$: 
\begin{equation}
 \sum_{k=1}^{n}A_{k}^{(\ell)}\to L
\text{ as $\ell\to \infty$}.
\label{eq:limit of summation}
\end{equation}
On the other hand, 
by \eqref{eq:arithmetic-mean inequality}, 
for $X,Y\in \mathcal{S}_{2}(H)_{++}$, we have 
%
%$$ 0\leq \left( (1-\lambda )X+\lambda Y+\frak{m}_{f}(X,Y)\right)^{1/2}
%\left((1-\lambda )X+\lambda Y-\frak{m}_{f}(X,Y)\right)
%\left( (1-\lambda )X+\lambda Y+\frak{m}_{f}(X,Y)\right)^{1/2}. $$
%%
%Then 
%
\begin{align*}
0 & \leq
{\rm Tr}\left[
\{(1-\lambda )X+\lambda Y+\frak{m}_{f}(X,Y)\}
\{(1-\lambda )X+\lambda Y-\frak{m}_{f}(X,Y)\}
\right]\\
& =
{\rm Tr}
\left\{(1-\lambda )X+\lambda Y
\right\}^{2}-{\rm Tr}\left(
\frak{m}_{f}(X,Y)^2\right)\\
& =
(1-\lambda){\rm Tr}(X^{2})+
\lambda {\rm Tr}(Y^{2})-
\lambda(1-\lambda){\rm Tr}
(X-Y)^{2}-
{\rm Tr}\left(
\frak{m}_{f}(X,Y)^2\right)\\
& =
(1-\lambda)\|X\|_{2}^{2}+
\lambda\|Y\|_{2}^{2}-
\lambda(1-\lambda)\|X-Y\|_{2}^{2}-\|\frak{m}_{f}(X,Y)\|_{2}^{2}.
\end{align*}
Hence 
$$ \|\frak{m}_{f}(X,Y)\|_{2}^{2}\leq
(1-\lambda)\|X\|_{2}^{2}+\lambda\|Y\|_{2}^{2}-
\lambda(1-\lambda)\|X-Y\|_{2}^{2}. $$
Therefore,
\begin{align*}
a_{\ell+1}& :=
\sum_{k=1}^{n}\| A_{k}^{(\ell+1)}\|_{2}^{2}\\
& \leq 
\sum_{k=1}^{n}\| A_{k}^{(\ell)}\|_{2}^{2}
-\lambda(1-\lambda)
\left(\sum_{k=1}^{n-1}\|A_{k}^{(\ell)}-
A_{k+1}^{(\ell)}\|_{2}^{2}+
\|A_{n}^{(\ell)}-
A_{1}^{(\ell)}\|_{2}^{2}\right)\\
& :=a_{\ell}-c_{\ell}\leq a_{\ell}.
\end{align*}
Since the numerical sequence $a_{\ell}$
is decreasing, it has a limit and it 
follows that $c_{\ell}\to 0$ as $\ell\to \infty$. Therefore, 
$$ \|A_{k}^{(\ell)}-A_{k+1}^{(\ell)}\|_{2}\to 0
\text{ and }
\|A_{n}^{(\ell)}-A_{1}^{(\ell)}\|_{2}\to 0 \text{ as $\ell\to \infty$}$$
for $k=1,2,...,n-1$. Combining formulas with 
\eqref{eq:limit of summation}, we have
$$ A_{k}^{(\ell)}\to \frac{1}{n}L 
\text{ as $\ell\to \infty$ for $k=1,2,...,n$.} $$
\end{proof}

The limit $L$ can be considered as 
a multivariate operator mean.
 
\begin{proof}[Proof of Theorem \ref{thm:matrices case}]
Assume that $T$ is a centered matrix.
By Lemma \ref{lem:centered matrices},
there exists a natural number $n$ such that
$|T|=U^{*n}|T|U^{n}$.
Let $A_{k}:=U^{*k}|T|U^{k}$ for 
$k=0,1,...,n-1$, and define 
$\{A_{k}^{(\ell)}\}$ by the same way 
to Theorem \ref{thm:multivariate operator mean}, i.e., 
$A_{k}^{(0)}:=A_{k}$ for $k=0,1,...,n-1$, 
$A_{k}^{(\ell+1)}:=
\frak{m}_{f}(A_{k+1}^{(\ell)}, A_{k}^{(\ell)})$
for $k=0,1,...,n-2$ and
$A_{n-1}^{(\ell+1)}:=
\frak{m}_{f}(A_{0}^{(\ell)}, A_{n-1}^{(\ell)})$.

Firstly, we shall prove that 
\begin{equation}
U^{*}A_{k}^{(\ell)}U=A_{k+1}^{(\ell)}, 
\text{ and }
U^{*}A_{n-1}^{(\ell)}U=A_{0}^{(\ell)}
\text{ for $k=0,1,...,n-2$}.
\label{eq: A_{k}}
\end{equation}
Since $A_{k}^{(0)}=A_{k}=U^{*k}|T|U^{k}$,
$$ U^{*}A_{k}^{(0)}U=U^{*(k+1)}|T|U^{k+1}
=A_{k+1}=A_{k+1}^{(0)}$$
holds for $k=0,1,...,n-2$, and
$U^{*}A_{n-1}^{(0)}U=A_{0}^{(0)}$
is shown similarly by $U^{*n}|T|U^{n}=|T|$.
Hence 
\eqref{eq: A_{k}} holds for $\ell=0$.
Assume that \eqref{eq: A_{k}} holds
for some $\ell$. Then 
\begin{align*}
U^{*}A_{k}^{(\ell+1)}U & =
U^{*} \frak{m}(A_{k+1}^{(\ell)}, 
A_{k}^{(\ell)})U\\
& =
 \frak{m}(U^{*}A_{k+1}^{(\ell)}U, 
U^{*}A_{k}^{(\ell)}U)\\
& =
 \frak{m}(A_{k+2}^{(\ell)}, A_{k+1}^{(\ell)})
 =
 A_{k+1}^{(\ell+1)},
 \end{align*}
and $U^{*}A_{n-1}^{(\ell+1)}U=A_{0}^{(\ell+1)}$ holds, too.
Hence we have \eqref{eq: A_{k}} 
for all $\ell=0,1,2,...$

Next, we shall prove that 
$|\Delta_{\frak{m}_{f}}^{\ell}(T)|=
A_{0}^{(\ell)}$ holds 
for $\ell=0,1,2,...$
by mathematical induction on $\ell$. 
The case $\ell =0$, we have
$$ |\Delta_{\frak{m}_{f}}^{0}(T)|=
|T|=A_{0}^{0}.$$
Assume that
$|\Delta_{\frak{m}_{f}}^{\ell}(T)|
=\frak{M}_{\ell}(U,|T|)=A_{0}^{(\ell)}$ for some $\ell$. Then 
\begin{align*}
|\Delta_{\frak{m}_{f}}^{\ell+1}(T)| & =
\frak{m}(U^{*}\frak{M}_{\ell}(U,|T|)U,
\frak{M}_{\ell}(U,|T|))\\
& =
\frak{m}(U^{*}A_{0}^{(\ell)}U, A_{0}^{(\ell)})\\
& =
\frak{m}(A_{1}^{(\ell)}, A_{0}^{(\ell)})
=
A_{0}^{(\ell+1)}.
\end{align*}
Hence $|\Delta_{\frak{m}_{f}}^{\ell}(T)|=
A_{0}^{(\ell)}$ for all $\ell=0,1,2,...$
By Theorem \ref{thm:multivariate operator mean}, $\{|\Delta_{\frak{m}_{f}}^{\ell}(T)|\}=\{A_{0}^{(\ell)}\}$ converges as $\ell\to \infty$.
By Theorem \ref{thm:iteration}, 
the polar decomposition of 
$\Delta^{\ell}_{\frak{m}_{f}}(T)$ is
$ \Delta^{\ell}_{\frak{m}_{f}}(T)=U|\Delta^{\ell}_{\frak{m}_{f}}(T)|$.
Therefore
$\{\Delta_{\frak{m}_{f}}^{\ell}(T)\}$ converges as $\ell\to \infty$.

We put $UL$ as the limit point of 
$\{\Delta_{\frak{m}_{f}}^{\ell}(T)\}$.
Then $\Delta_{\frak{m}_{f}}(UL)=UL$ holds. Hence 
$UL$ is normal by \cite[Theorem 2.4]{Y 2021}.
\end{proof}

%%%%%%%%%%%%%%%%%%%%%%%%%%%%%%%%%%%%%%%%%%%%%%%%%%%
\section{Standard form and its applications}

Let $M$ be a von Nuemann algebra and $\phi$ be a weight of $M$. In the case that $M$ is finite, we may consider $\phi$ as a faithful state. In our case we consider $M = \mathcal{M}_m(\C)$. So, $\phi$ be a normalized trace on $\mathcal{M}_m(\C)$, that is, $\phi = \frac{1}{m}\Tr$.

Let $(H_\phi, \pi_\phi, \eta_\phi)$ be the GNS construction of $M$. Set $N_\phi = \{a \in M| \phi(a^*a) < \infty\}$. Then, $N_\phi$ is $\sigma$-weak dense in $M$ and $\eta_\phi\colon N_\phi \rightarrow H_\phi$ is a linear map such that $\overline{\eta_\phi(N_\phi)} = H_\phi$. Then, $\pi_\phi(a)\eta_\phi(b) = \eta_\phi(ab)$. Then, $\pi_\phi\colon M \rightarrow B(H_\phi)$ is a normal representation. Set  a conjugate-linear operator $S_\phi$ on  a dense subspace $\eta_\phi(N_\phi \cap N_\phi^*)$ of $H_\phi$ by $S_\phi(\eta_\phi(a)) = \eta_\phi(a^*)$, which is closable.
By taking the polar decomposition $\overline{S_\phi} = J_\phi\Delta^{1/2}_\phi$ of $\overline{S_\phi}$ which is the closure of $S_\phi$, we define a positive self-adjoint operator $\Delta_\phi$ called the modular operator, and a conjugate involution $J_\phi$ ($J_\phi^2 = 1$) called the modular conjugation. Then by the Tomita's fundamental theorem \cite{Bl 2006} says that

\begin{enumerate} 
\item
$J_\phi M J_\phi = M'$,
\item
$\Delta^{it}_\phi M \Delta^{-it}_\phi = M$ $(t \in \R)$.
\end{enumerate}

\vskip 2mm 

In the case of $M = \mathcal{M}_m(\C)$. Then, $H_\phi = \mathcal{M}_m(\C)$ as a set and $S_\phi = J_\phi$. We do not need the general theory when $M = \mathcal{M}_{m}(\C)$.

\vskip 2mm

\begin{lem}\label{lem:spacial automorphism}
Let $(H_\phi, \pi_\phi, \eta_\phi)$ be the standard representation of $\mathcal{M}_m(\C)$, $U \in \mathcal{M}_m(\C)$ a unitary. Define $\alpha_U\in \mathcal{B}(H_{\phi})$ by $\alpha_U(\eta_\phi(a)) = \eta_\phi(U^*aU)$. Then, $(\alpha_U)^* = \alpha_{U^*}$ and  $\|\alpha_U\| = 1$. 
\end{lem}

\vskip 2mm

\begin{proof}
For any $a, b \in \mathcal{M}_m(\C)$
\begin{align*}
(\alpha_U(\eta_\phi(a))|\eta_\phi(b)) &= (\eta_\phi(U^*aU)|\eta_\phi(b))\\
&= \phi(b^*U^*aU) = \phi(Ub^*U^*a)\\
&= (\eta_\phi(a)|\alpha_{U^*}(\eta_\phi(b)))
\end{align*}

\vskip 2mm

Hence, $(\alpha_U)^* = \alpha_{U^*}$.

For $a \in \mathcal{M}_m(\C)$
$$
\|\alpha_U(\eta_\phi(a))\|^2 = \|\eta_\phi(U^*aU)\|^2 = \phi(U^*a^*aU)
\leq \|a\|^2.
$$
Hence, $\|\alpha_U\| \leq 1$. Since $\alpha_U(1) = 1$, $\|\alpha_U\| = 1$.
\end{proof}

\vskip 2mm

\begin{lem}
Let $U \in \mathcal{M}_m(\C)$ be a unitary. Then, $\ker (\alpha_U - I) = \{U\}^{'}$, the relative commutant of $U$.
\end{lem}

\vskip 2mm

\begin{proof}
For any $a \in \ker(\alpha_U- I)$. Then, $(\alpha_U-I)(a) = 0$ . That is, $U^*aU = a$. Hence, $aU = Ua$, $a \in \{U\}^{'}$. The converse inclusion is easy.
\end{proof}

\vskip 2mm

The following is the natural extension of \cite[Theorem~4.1]{DS 2009}.

\vskip 2mm

\begin{prp}\label{prp:extension of Dykema}
Let $T \in \mathcal{B}(H)$ be a contraction, i.e., $\|T\|\leq 1$ and $\lambda \in (0, 1)$.  
Then, for every $v \in H$, the vector 
$$
L_n(\lambda)(v) = \sum_{k=0}^n 
\left(\begin{array}{c}
n\\
k
\end{array}
\right) 
(1 - \lambda)^k\lambda^{n-k}T^kv
$$
converges to $Pv$ as $n \rightarrow \infty$, where $P$ is the orthogonal projection of $H$ onto its subspace $\ker (T -I) = \{x \in H | Tx = x\}$.
\end{prp}
The case $\lambda =\frac{1}{2}$ has been shown in  \cite[Theorem~4.1]{DS 2009}. Hence Proposition \ref{prp:extension of Dykema} is a generalization of \cite[Theorem~4.1]{DS 2009}. However, the proof is simpler than \cite[Theorem~4.1]{DS 2009}.

\begin{proof}
Set $a_{n,k}:=(1-\lambda)^{k}\lambda^{n-k}$ $(0\leq k\leq n)$. Note that
$ L_{n}(\frac{1}{2})=H_n$ in \cite[Theorem~4.1]{DS 2009} and $\sum_{k=0}^{n}\begin{pmatrix} n \\ k \end{pmatrix}a_{n,k}=1$.

As in the same argument of \cite[Theorem~4.1]{DS 2009}, we may assume that $v=Ty-y$ for some $y\in H$. 
\begin{align*}
L_{n}(\lambda)(v) & =
\sum_{k=0}^{n}\begin{pmatrix} n \\ k \end{pmatrix} a_{n,k}(T^{k+1}y-T^{k}y)\\
& =
-a_{n,0}y+\sum_{k=1}^{n}\left[
\begin{pmatrix} n \\ k-1 \end{pmatrix}a_{n,k-1}-
\begin{pmatrix} n \\ k \end{pmatrix}a_{n,k}\right]T^{k}y+a_{n,n}T^{n+1}y.
\end{align*}
Since $| \begin{pmatrix} n\\ k-1\end{pmatrix}a_{n,k-1}-\begin{pmatrix} n\\ k \end{pmatrix}a_{n,k}|=\frac{|k-(1-\lambda)(n+1)|}{k\lambda} \begin{pmatrix} n \\ k-1 \end{pmatrix} a_{n,k-1}$, we get
\begin{align*}
\| L_{n}(\lambda)(v)\| &\leq
\left( a_{n,0}+\sum_{k=1}^{n}\frac{|k-(1-\lambda)(n+1)|}{k\lambda} \begin{pmatrix} n \\ k-1 \end{pmatrix} a_{n,k-1}+a_{n,n}\right)\|y\| \\
& =
\left(\frac{1-\lambda}{\lambda} n a_{n,0}+
\sum_{k=1}^{n} \frac{|(k+1)-(1-\lambda)(n+1)|}{\lambda(k+1)} \begin{pmatrix} n \\ k \end{pmatrix} a_{n,k}\right) \|y\|, 
\end{align*}
where the last equation holds for $n\geq \frac{\lambda}{1-\lambda}$, i.e., it holds for sufficiently large $n$.

Since $\lim_{n\to \infty}\frac{1-\lambda}{\lambda} n a_{n,0}=\lim_{n\to \infty}n(1-\lambda)\lambda^{n-1}=0$, we may show that
$$ \lim_{n\to \infty} \sum_{k=1}^{n}\frac{|(k+1)-(1-\lambda)(n+1)|}{\lambda(k+1)} \begin{pmatrix} n \\ k \end{pmatrix} a_{n,k}=0.
$$

Put $M=\lfloor (1-\lambda)n-\lambda\rfloor$. Then 
\begin{align*}
& \sum_{k=1}^{n} \frac{|(k+1)-(1-\lambda)(n+1)|}{\lambda(k+1)} \binom{n}{k} a_{n,k}\\
& =
\sum_{k=1}^{M} \frac{(1-\lambda)(n+1)-(k+1)}{\lambda(k+1)} \binom{n}{k} a_{n,k} +
\sum_{k=M+1}^{n} \frac{(k+1)-(1-\lambda)(n+1)}{\lambda(k+1)} \binom{n}{k} a_{n,k}\\
& =
\sum_{k=1}^{M} \left[
\frac{(n+1)!}{(k+1)!(n-k)!}(1-\lambda)^{k+1}\lambda^{n-k-1} -
\frac{n!}{k!(n-k)!}(1-\lambda)^{k}\lambda^{n-k-1}
\right] \\
& \qquad +
\sum_{k=M+1}^{n} \left[
\frac{n!}{k!(n-k)!}(1-\lambda)^{k}\lambda^{n-k-1} -
\frac{(n+1)!}{(k+1)!(n-k)!}(1-\lambda)^{k+1}\lambda^{n-k-1}
\right]\\
& =
\sum_{k=1}^{M} \left[
\binom{n+1}{k+1} (1-\lambda)^{k+1}\lambda^{n-k-1} -
\binom{n}{k} (1-\lambda)^{k}\lambda^{n-k-1}
\right]\\
& \qquad +
\sum_{k=M+1}^{n} \left[
\binom{n}{k} (1-\lambda)^{k}\lambda^{n-k-1} -
\binom{n+1}{k+1} (1-\lambda)^{k+1}\lambda^{n-k-1}
\right]\\
& =
\sum_{k=1}^{M} \left[
\binom{n}{k} (1-\lambda)^{k+1}\lambda^{n-k-1} +
\binom{n}{k+1} (1-\lambda)^{k+1}\lambda^{n-k-1} -
\binom{n}{k} (1-\lambda)^{k}\lambda^{n-k-1}
\right]\\
& \qquad +
\sum_{k=M+1}^{n} \left[
\binom{n}{k} (1-\lambda)^{k}\lambda^{n-k-1} -
\binom{n}{k} (1-\lambda)^{k+1}\lambda^{n-k-1} -
\binom{n}{k+1} (1-\lambda)^{k+1}\lambda^{n-k-1}
\right]\\
& =
\sum_{k=1}^{M} \left[
\binom{n}{k+1} (1-\lambda)^{k+1}\lambda^{n-k-1} -
\binom{n}{k} (1-\lambda)^{k}\lambda^{n-k}
\right]\\
& \qquad +
\sum_{k=M+1}^{n} \left[
\binom{n}{k} (1-\lambda)^{k}\lambda^{n-k} -
\binom{n}{k+1} (1-\lambda)^{k+1}\lambda^{n-k-1}
\right]\\
& =
\binom{n}{M+1} (1-\lambda)^{M+1}\lambda^{n-M-1} -
\binom{n}{1} (1-\lambda)\lambda^{n-1}\\
& \qquad +
\binom{n}{M+1} (1-\lambda)^{M+1}\lambda^{n-M-1} -
\binom{n}{n+1} (1-\lambda)^{n+1}\lambda^{-1}\\
& =
2 \binom{n}{M+1} (1-\lambda)^{M+1}\lambda^{n-M-1} -
n (1-\lambda) \lambda^{n-1}.
\end{align*}
where $\begin{pmatrix} n \\ n+1 \end{pmatrix}=0$. Since $\lim_{n\to \infty} n(1-\lambda)\lambda^{n-1}=0$, we only prove
$$ \lim_{n\to \infty} \begin{pmatrix} n \\ M+1 \end{pmatrix}
 (1-\lambda)^{M+1}\lambda^{n-M-1}=0$$
 for $M=\lfloor (1-\lambda)n-\lambda\rfloor$.
 
 We note that since $(1-\lambda)n-\lambda-1\leq M\leq (1-\lambda)n-\lambda$,
 we have $\lim_{n\to \infty}\frac{M}{n}=1-\lambda$. By Starling's formula, 
$$ \sqrt{2\pi}n^{n+1/2}e^{-n}\leq n!\leq
e n^{n+1/2} e^{-n}$$
holds for all natural number $n$. Then we have
\begin{align*}
0 & \leq 
\begin{pmatrix} n \\ M+1 \end{pmatrix}
(1-\lambda)^{M+1}\lambda^{n-M-1}
= 
\frac{n!}{(M+1)!(n-M-1)!} (1-\lambda)^{M+1}\lambda^{n-M-1}\\
& \leq
\frac{e n^{n+1/2}e^{-n}}{\sqrt{2\pi}(M+1)^{M+1+1/2}e^{-M-1}\sqrt{2\pi}(n-M-1)^{n-M-1+1/2}e^{-n+M+1}}
(1-\lambda)^{M+1}\lambda^{n-M-1}\\
& =
\frac{e}{2\pi}
\frac{n^{n}}{(M+1)^{M+1}(n-M-1)^{n-M-1}}
\frac{\sqrt{n}}{\sqrt{M+1}\sqrt{n-M-1}}
(1-\lambda)^{M+1}\lambda^{n-M-1}\\
& =
\frac{e}{2\pi}\frac{1}{\sqrt{n}}
\frac{1}{\sqrt{(\frac{M}{n}+\frac{1}{n})(1-\frac{M}{n}-\frac{1}{n})}}
\left(\frac{(1-\lambda)n}{M+1}\right)^{M+1}
\left(\frac{\lambda n}{n-M-1}\right)^{n-M-1}\\
& \leq
\frac{e}{2\pi}\frac{1}{\sqrt{n}}
\frac{1}{\sqrt{(\frac{M}{n}+\frac{1}{n})(1-\frac{M}{n}-\frac{1}{n})}}
\left(\frac{(1-\lambda)n}{(1-\lambda)n-\lambda}\right)^{M+1}
\left(\frac{\lambda n}{n-(1-\lambda)n+\lambda-1}\right)^{n-M-1}\\
& \hspace*{5cm} \text{(by $(1-\lambda)n-\lambda-1\leq M\leq (1-\lambda)n-\lambda$)}\\
& =
\frac{e}{2\pi}\frac{1}{\sqrt{n}}
\frac{1}{\sqrt{(\frac{M}{n}+\frac{1}{n})(1-\frac{M}{n}-\frac{1}{n})}}
\left(1-\frac{\lambda}{(1-\lambda)n}\right)^{-n(\frac{M}{n}+\frac{1}{n})}
\left(1-\frac{1-\lambda}{\lambda n}\right)^{-n(1-\frac{M}{n}-\frac{1}{n})}\\
& \to
\frac{e^2}{2\pi}\cdot 0\cdot \frac{1}{\sqrt{\lambda(1-\lambda)}}=0
\end{align*}
as $n\to \infty$. 

Hence the proof is completed.
\end{proof}

%\begin{proof}
%Set $a_{n,k} = (1 - \lambda)^k\lambda^{n-k}$ $(0 \leq k \leq n)$.
%Note that $L_n(\frac{1}{2}) = H_n$ in \cite[Theorem~4.1]{DS 2009} and $\sum_{k=0}^n 
%\left(\begin{array}{c}
%n\\
%k
%\end{array}
%\right) 
%a_{n,k}
%= 1$.
%

\vskip 2mm

\begin{cor}\label{cor:extension of Dykema}
Let $T \in \mathcal{B}(H)$ be a contraction, i.e., $\|T\|\leq 1$ and $\lambda \in (0, 1)$. %$0 < 1- \lambda \leq \frac{1}{2}$. 
%Set $a_k = \lambda^k(1 - \lambda)^{n-k}$ $(0 \leq k \leq n)$.
Then, for every $v \in H$, there exists $\lim_{n\rightarrow\infty}L_n(\lambda)(v)$, and it does not depend on $\lambda$.
\end{cor}

\begin{proof}
%For every $v \in H$ and $0 \leq \lambda \leq 1$, 
It follows from Proposition~\ref{prp:extension of Dykema}.
%that $\lim_{n\rightarrow\infty}L_n(\lambda)(v) = Pv$, where $P\colon H \rightarrow \{x \in H| Tx = x\}$ is a projection.
\end{proof}

\vskip 2mm

\vskip 2mm

\begin{thm}\label{thm:weighted arithmetic}
Let $H$ be a finite dimensional Hilbert space and $\lambda \in (0, 1)$, 
$f_\lambda$ be a function by $f_\lambda(x) = 1 - \lambda + \lambda x$. Then, for any invertible element $T \in \mathcal{B}(H)$ with the polar decomposition $T = U|T|$,
$\lim_{n\rightarrow\infty}\Delta_{f_{\lambda}}^n(T)$ exists, and it does not depend on $\lambda$.
\end{thm}

\vskip 2mm

\begin{proof}
Let $(H_\phi, \pi_\phi, \eta_\phi)$ be the standard form of $\mathcal{M}_m(\C)$ and $\alpha_U \in \mathcal{B}(H_\phi)$ as in Lemma~\ref{lem:spacial automorphism}. Then,
by Proposition \ref{prp: power mean_centered}.

\begin{align*}
\eta_\phi(\Delta_{f_{\lambda}}^n(T)) & = \eta_\phi\left(U\left[\sum_{k=0}^n
\left(\begin{array}{c}
n\\
k
\end{array}
\right)
a_{n,k}
U^{*k}|T|U^k \right]
\right) \\
& =
\pi_\phi(U)\left(\sum_{k=0}^n
\left(\begin{array}{c}
n\\
k
\end{array}
\right)
a_{n,k}
\alpha_U^k
\right)
(\eta_\phi(|T|)).
\end{align*}

Note that $\displaystyle \lim_{n\rightarrow\infty}\left(\sum_{k=0}^n
\left(\begin{array}{c}
n\\
k
\end{array}
\right)
a_{n,k}
\alpha_U^k
\right)
(\eta_\phi(|T|)) = P\eta_\phi(|T|)$ for a projection $P\colon H_\phi \rightarrow \eta_\phi(\{U\}^{'})$ by Proposition~\ref{prp:extension of Dykema}, where $\{U\}^{'} = \{X \in \mathcal{M}_m(\C)| UX = XU\}$.

Since the norm $\|\quad\|_\phi$ is equivalent to the standard norm $\| \quad\|$ on $\mathcal{M}_m(\C)$, $\lim_{n\rightarrow\infty}\Delta_{f_{\lambda}}^n(T)$ exists, and $\lim_{n\rightarrow\infty}\Delta_{f_{\lambda}}^n(T) =UP$.
Moreover we obtain that 
the limit point does not depend on $\lambda$
by Corollary \ref{cor:extension of Dykema}.
\end{proof}

\vskip 2mm

The following corollary is a generalization
of \cite[Theorem~5.1]{Y 2021}.

\begin{cor}\label{cor:iteration of arithmetric mean}
Let $\lambda \in (0, 1)$, and $f_{\lambda}$ be a function by $f_{\lambda}(x)=1-\lambda +\lambda x$. For an invertible $T\in \mathcal{M}(\mathbb{C})$, let $T=U|T|$
be the polar decomposition. Then, there exists
$\lim_{n\rightarrow\infty}\Delta_{f_\lambda}^n(T)$, 
and it does not depend on $\lambda$.
\end{cor}

\begin{proof}
It follows from the proof in Theorem~\ref{thm:weighted arithmetic} and Corollary \ref{cor:extension of Dykema}.
\end{proof}

\vskip 2mm

\begin{rmk}\label{rrnk:iteration in finite rank}
Let $(H_\phi, \pi_\phi, \eta_\phi)$ be the standard form of $\mathcal{M}_m(\C)$ and $U \in B(H_\phi)$ be a unitary.
Then, for any $S \in \mathcal{M}_m(\C)$,  
$\displaystyle \left\{\sum_{k=0}^n
\left(\begin{array}{c}
n\\
k
\end{array}
\right)
a_{n,k}
U^{*k}SU^k
\right\}_{n\in\N}
$ 
converges as $n\to \infty$
by the same argument in Theorem~\ref{thm:weighted arithmetic}.
\end{rmk}

\vskip 2mm

\begin{cor}\label{cor:matrix harmonic mean}
Let $H$ be a finite dimensional Hilbert space,
$\lambda \in (0,1)$, $r\in [-1,1)\setminus\{0\}$ and define
$f_{r,\lambda}(x)=[1-\lambda+\lambda x^{r}]^{1/r}$ 
on $x>0$.
%$0 < 1- \lambda \leq \frac{1}{2}$. 
Let $T = U|T|$ be the polar decomposition of an invertible centered matrix $T \in \mathcal{M}_m(\C)$. 
Then, $\lim_{n\rightarrow\infty}
\Delta_{f_r,\lambda}^n(T)$ exists, and 
 it does not depend on $\lambda$.
\end{cor}

\begin{proof}
By Proposition \ref{prp: power mean_centered},
we have the polar decomposition of 
$\Delta_{f_{r,\lambda}}^n(T)$ as follows:
$$ \Delta_{f_{r,\lambda}}^n(T)=
U \left[ \sum_{k=0}^{n}
\begin{pmatrix} n \\ k \end{pmatrix}
(1-\lambda)^{k}\lambda^{n-k} 
U^{*k}|T|^{r}U^{k}\right]^{1/r}. $$
By Remark \ref{rrnk:iteration in finite rank}, $\{|\Delta_{f_{r,\lambda}}^n(T)|^{r}\}$ converges as $n\to \infty$.
Moreover, by Corollary \ref{cor:iteration of arithmetric mean},
the limit point does not depend on $\lambda$.
%
%$$ \lim_{n\to\infty}|\Delta_{f_{r,\lambda}}^n(T)|^{r}=
%\lim_{n\to\infty}|\Delta_{f_{r,1/2}}^n(T)|^{r}. $$
%
The proof is completed.
\end{proof}

%%%%%%%%%%%%%%%%%%%%%%%%%%%%%%%%%%%%%%%%%%%%%%%%%%%%%%%%
\section{Compact operators case}
In this section, we shall show that IAS 
converges to a normal operator in the case of centered compact operators.

\begin{prp}\label{prp:finite rank}
Let $H$ be an infinite dimensional Hilbert space, $U \in \mathcal{B}(H)$ be a unitary %(or an isometry) 
and $\lambda \in (0, 1)$. %$0 < 1 - \lambda \leq \frac{1}{2}$ 
%and $a_k = (1-\lambda)^k\lambda^{n-k}$ $(0 \leq k \leq n)$. 
Then, for any $F \in \mathcal{F}(H)_{+}$, 
there exists
$$
\lim_{n\rightarrow\infty}\sum_{k=1}^n
\left(\begin{array}{c}
n\\
k
\end{array}
\right)
(1-\lambda)^k\lambda^{n-k}
U^{*k}FU^k 
 \in \mathcal{F}(H),$$ 
and it does not depend on $\lambda$.
\end{prp}

\begin{proof}
Set $a_{n,k} = (1-\lambda)^k\lambda^{n-k}$.
Let $\|\quad\|_2$ be the Schatten 2-norm. Then for any $F \in \mathcal{F}(H)_+$, $\|F\| \leq \|F\|_2$. Hence, from the same argument in Theorem~\ref{thm:weighted arithmetic} (here $H_\phi = \mathcal{S}_2(H)$) 
$$
\lim_{n\rightarrow\infty}
\sum_{k=0}^n
\left(\begin{array}{c}
n\\
k\\
\end{array}
\right)
a_{n,k}
U^{*k}FU^k$$
exists. Moreover, from Proposition \ref{prp:extension of Dykema} 
%\cite[Theorem~4.1]{DS 2009} 
there exists a projection $P \in \mathcal{B}(\mathcal{S}_2(H))$ such that 
$\displaystyle 
\lim_{n\rightarrow\infty}
\sum_{k=0}^n
\left(\begin{array}{c}
n\\
k\\
\end{array}
\right)
a_{n,k}
U^{*k}FU^k
= PF \in \mathcal{F}(H)$.
\end{proof}

\begin{cor}\label{cor:compact}
Let $H$ be an infinite dimensional Hilbert space, $U\in \mathcal{B}(H)$ be a unitary, $\lambda \in (0, 1)$ %$0 < 1 - \lambda \leq \frac{1}{2}$, 
%$a_k = (1-\lambda)^k\lambda^{n-k}$ $(0 \leq k \leq n)$ 
and $K\in \mathcal{K}(H)_{+}$. Then, 
there exists
$$
\lim_{n\rightarrow\infty}\sum_{k=0}^n
\left(\begin{array}{c}
n\\
k\\
\end{array}
\right)
(1-\lambda)^k\lambda^{n-k}
U^{*k}KU^k 
$$
and it does not depend on $\lambda$.
\end{cor}

\begin{proof}
Let $a_{n,k}=(1-\lambda)^k \lambda^{n-k}$.
For any $\epsilon > 0$ there exists $F \in \mathcal{F}(H)_{+}$ such that $\|K - F\| < \epsilon$.

Since $\displaystyle \left\{\sum_{k=0}^n 
\left(\begin{array}{c}
n\\
k
\end{array}
\right)
a_{n,k}
U^{*k}FU^k
\right\}_{n\in\N}
$ 
has a limit as $n\to \infty$
by Proposition~\ref{prp:finite rank}, there exists $n_0 \in \N$ such that for any $n, m \geq n_0$ 
$$
\|\sum_{k=0}^n
\left(\begin{array}{c}
n\\
k
\end{array}
\right)
a_{n,k}
U^{*k}FU^k
-
\sum_{k=0}^m
\left(\begin{array}{c}
m\\
k
\end{array}
\right)
a_{m,k}
U^{*k}FU^k
\| < \epsilon.
$$
Then, we have 

\begin{align*}
&\|\sum_{k=0}^n
\left(\begin{array}{c}
n\\
k
\end{array}
\right)
a_{n,k}
U^{*k}KU^k
-
\sum_{k=0}^m
\left(\begin{array}{c}
m\\
k
\end{array}
\right)
a_{m,k}
U^{*k}KU^k\|\\
%&=
%\|\frac{1}{2^n}\sum_{k=0}^n
%\left(\begin{array}{c}
%n\\
%k
%\end{array}
%\right)
%U^{*k}KU^k
%-
%\sum_{k=0}^m
%\left(\begin{array}{c}
%m\\
%k
%\end{array}
%\right)
%a_k
%U^{*k}K U^k\|\\
&= 
\|\sum_{k=0}^n
\left(\begin{array}{c}
n\\
k
\end{array}
\right)
a_{n,k}
U^{*k}(K- F)U^k\\
& \hspace*{2cm}
+
\sum_{k=0}^n
\left(\begin{array}{c}
n\\
k
\end{array}
\right)
a_{n,k}
U^{*k}FU^k
-
\sum_{k=0}^m
\left(\begin{array}{c}
m\\
k
\end{array}
\right)
a_{m,k}
U^{*k}F U^k\\
& \hspace*{3cm} -
\sum_{k=0}^m
\left(\begin{array}{c}
m\\
k
\end{array}
\right)
a_{m,k}
U^{*k}(K-F) U^k\|\\
&\leq 2\|K - F\| + \|\sum_{k=0}^n
\left(\begin{array}{c}
n\\
k
\end{array}
\right)a_{n,k}
U^{*k}FU^k
-
\sum_{k=0}^m
\left(\begin{array}{c}
m\\
k
\end{array}
\right)
a_{m,k}
U^{*k}F U^k\|\\
&< 3 \epsilon
\end{align*}
This implies that 
$\displaystyle \lim_{n\rightarrow\infty}\sum_{k=0}^n
\left(\begin{array}{c}
n\\
k\\
\end{array}
\right)
a_{n,k}
U^{*k}KU^k \in \mathcal{K}(H)_{+}$ 
exists because the sequence is 
a Cauchy sequence.

By Proposition~\ref{prp:finite rank}, we get the conclusion.
\end{proof}

\vskip 2mm

\vskip 2mm

Using Proposition \ref{prp: power mean_centered}
and Corollary~\ref{cor:compact},
we have the following result.

\begin{thm}\label{thm:convergence power mean}
For $\lambda\in (0,1)$ and $r\in [-1,1]\setminus\{0\}$, let
$f_{r,\lambda}(x)=
[1-\lambda+\lambda x^{r}]^{1/r}$.
Assume that $T=U|T| \in \mathcal{B}(H)$
is invertible and centered.
If $T\in \mathcal{K}(H)$ or a $n$--by--$n$ 
matrix, then there exists $\lim_{n\to \infty}\Delta_{f_{r,\lambda}}^{n}(T)$.
Moreover, the limit point 
does not depend on $\lambda$.
\end{thm}

\begin{proof}
By Proposition \ref{prp: power mean_centered},
we have the polar decomposition of 
$\Delta_{f_{r,\lambda}}^n(T)$ as follows:
$$ \Delta_{f_{r,\lambda}}^n(T)=
U \left[ \sum_{k=0}^{n}
\begin{pmatrix} n \\ k \end{pmatrix}
(1-\lambda)^{k}\lambda^{n-k} 
U^{*k}|T|^{r}U^{k}\right]^{1/r}. $$
By Corollary~\ref{cor:compact}, $\{|\Delta_{f_{r,\lambda}}^n(T)|^{r}\}$ converges as $n\to \infty$.
Moreover, by Corollary~\ref{cor:compact},
the limit point does not depend on $\lambda$.
%
%$$ \lim_{n\to\infty}|\Delta_{f_{r,\lambda}}^n(T)|^{r}=
%\lim_{n\to\infty}|\Delta_{f_{r,1/2}}^n(T)|^{r}. $$
%
The proof is completed.
\end{proof}

\vskip 2mm

\begin{rmk}
If $r=1$ in Theorem \ref{thm:convergence power mean}, then the centered condition of $T$ is not needed.
See Remark \ref{rm:r=1 case}.
\end{rmk}

\vskip 2mm
%%%%%%%%%%%%%%%%%%%%%%%%%%%%%%%%%%%%%%%%%%%%%%%%%%%%%%
\section{$\mathcal{AN}$ and $\mathcal{AM}$ operators}

In this section, we shall consider larger classes of operators than the class of compact operators.
We shall treat $\mathcal{AN}$ and $\mathcal{AM}$ operators.
Let $H_1$ and $H_2$ be infinite dimensional Hilbert spaces. We recall that $T \in \mathcal{B}(H_1, H_2)$ is said to be an {\it absolutely norm attaining operator} or an $\mathcal{AN}$ operator for short, if $T|_M$, the restriction of $T$ to $M$, is norm attaining for every non-zero closed subspace $M\subseteq H_1$, that is, there exists a unit vector $x \in M$ such that $\|T|_M\| = \|T|_M x\|$
\cite{CN 2012}.
We denote $\mathcal{AN}(H_1,H_2)$ as the class of 
$\mathcal{AN}$ operators in $\mathcal{B}(H_1,H_2)$. Especially,
if $H_1=H_2$, then we write $\mathcal{AN}(H)$ for short.
An operator $T \in \mathcal{B}(H_1, H_2)$ is said to be an {\it absolutely minimal attaining operator} or an $\mathcal{AM}$ operator for short, if $T|_M$, the restriction of $T$ to $M$, is minimal attaining for every non-zero closed space $M\subseteq H_1$, that is, there exists a unit vector $x \in M$ such that $m(T|_M) = \|T|_M x\|$, where $m(T) = \inf\{\|Tx\| \ |\  x \in S_{H_1}\}$ \cite{CN 2014} and \cite[Definition 2.4]{BR2020}.
We denote $\mathcal{AM}(H_1, H_2)$ as the class of $\mathcal{AM}$ operators in $\mathcal{B}(H_1,H_2)$. Especially, if $H_1=H_2$, then we write $\mathcal{AM}(H)$ for short.
The class of $\mathcal{AN}$ operators includes all compact operators and
partial isometries with finite dimensional kernels.
The class of $\mathcal{AM}$ operators includes finite rank operators and partial isometries with finite dimensional kernels or ranges. We note that the class of $\mathcal{AM}$ operators includes non-injective compact operators.
Essential spectrum of every 
$\mathcal{AN}$ or $\mathcal{AM}$ operator is a singleton
\cite[Theorem 2.4]{GR 2018} and \cite[Theorem 3.10]{BR2020}, respectively.
Properties of $\mathcal{AN}$ and $\mathcal{AM}$ operators are introduced in \cite{BR2020, CN 2012, CN 2014, GRS2018, GOUY pre2, NR 2019, PP 2017, Ramesh 2014, GR 2018, GS 2021, GOUY 2022} and reference there in.

In this section, we shall show that the IASs of invertible $\mathcal{AN}$ and $\mathcal{AM}$ operators with respect to the power means converge. Moreover, the limit operation preserves $\mathcal{AN}$ and $\mathcal{AM}$ properties, respectively. 

We shall start this section by introducing the following basic properties of $\mathcal{AN}$ and $\mathcal{AM}$ operators.
The following results are shown in 
\cite[Theorem 2.5]{NR 2019}, \cite[Theorem 5.9]{GRS2018} and \cite[Theorem 5.6]{GRS2018}.

%\vskip 2mm

\begin{thm}\label{basic characterizations}
The following statements hold.
\begin{itemize}
\item[(i)] $T\in \mathcal{AN}(H)_+$ if and only if there exist a positive number $\alpha\geq 0$, $K\in \mathcal{K}(H)_+$ and a self-adjoint $F\in \mathcal{F}(H)$ such that
$$ T=\alpha I+K+F \quad \mbox{(\cite[Theorem 2.5]{NR 2019}).}$$
\item[(ii)] $T\in \mathcal{AM}(H)_+$ if and only if there exist a positive number $\alpha\geq 0$, $K\in \mathcal{K}(H)_+$ and $F\in \mathcal{F}(H)_+$ such that
$$ T=\alpha I-K+F, $$
$KF=0$ and $0\leq K\leq \alpha I$. The decomposition is unique (\cite[Theorem 5.9]{GRS2018}).
\item[(iii)] Let $F\in \mathcal{F}(H)$ be self-adjoint and $K\in \mathcal{K}(H)_+$. Then for every positive number $\alpha$ satisfying $\frac{\|K\|}{2}\leq \alpha$, we have 
$$\alpha I-K+F\in \mathcal{AM}(H)\quad \mbox{(\cite[Theorem 5.6]{GRS2018}).} $$
\end{itemize}
\end{thm}

\begin{lem}\label{lem:invertibility}
Let $T \in B(H)$ be invertible. Then the following statements hold.
\begin{itemize}
\item[(i)] $T \in \mathcal{AN}(H)$ if and only if  $T^{-1} \in \mathcal{AM}(H)$,
\item[(ii)] $T \in \mathcal{AN}(H)$ if and only if $T^*\in \mathcal{AN}(H)$.
\end{itemize}
\end{lem}

\begin{rmk}
In general, $T\in \mathcal{AN}(H)$ does not imply $T^*\in \mathcal{AN}(H)$ \cite[Proposition 3.14]{CN 2012}.
\end{rmk}

\begin{proof}[Proof of Lemma \ref{lem:invertibility}]
(i) Let $T \in \mathcal{AN}(H)$. Then 
\begin{align*}
T \in \mathcal{AN}(H)
&\leftrightarrow |T| \in \mathcal{AN}(H) \quad  (\hbox{\cite[Lemma~6.2]{PP 2017}})\\
&\leftrightarrow |T|^{-1} \in \mathcal{AM}(H) \quad  (\hbox{\cite[Theorem~3.8]{BR2020}})\\
&\leftrightarrow |(T^*)^{-1}| \in \mathcal{AM}(H) \quad (|T|^{-1} = |(T^*)^{-1}|)\\
&\leftrightarrow (T^*)^{-1} \in \mathcal{AM}(H) \quad (\hbox{\cite[Proposition~3.2]{GRS2018}})\\
&\leftrightarrow T^{-1} \in \mathcal{AM}(H) \quad (\hbox{\cite[Corollary~3.15]{BR2020}}).
\end{align*}

(ii) Suppose that $T \in \mathcal{AN}(H)$. %Then, 
By (i), it is equivalent to $T^{-1} \in \mathcal{AM}(H)$.
Hence, 
\begin{align*}
T^{-1} \in \mathcal{AM}(H)
&\rightarrow |T^{-1}| \in \mathcal{AM}(H) \quad (\hbox{\cite[Proposition~3.2]{NR 2019}})\\
&\rightarrow |T^*|^{-1} \in \mathcal{AM}(H)  \quad (|T|^{-1} = |(T^*)^{-1}|)\\
&\rightarrow |T^*| \in \mathcal{AN}(H) \quad (\hbox{\cite[Corollary~3.15]{BR2020}})\\
&\rightarrow T^* \in \mathcal{AN}(H) \quad  (\hbox{\cite[Lemma~6.2]{PP 2017}})
\end{align*}
\end{proof}

%
%Conversely,
%
%\begin{align*}
%T^{-1} \in \mathcal{AM}(H) &\rightarrow (T^*)^{-1} \in \mathcal{AM}(H) \quad  (\hbox{\cite[Corollary~3.15]{BR2020}})\\
%&\rightarrow |(T^*)^{-1}| \in \mathcal{AM}(H) \quad  (\hbox{\cite[Proposition~3.2]{GRS2018}})\\
%&\rightarrow  |T|^{-1} \in \mathcal{AM}(H) \quad (|T|^{-1} = |(T^*)^{-1}|) \\
%&\rightarrow |T| \in \mathcal{AN}(H) \quad  (\hbox{\cite[Theorem~3.8]{BR2020}})\\
%&\rightarrow T \in \mathcal{AN}(H) \quad (\hbox{\cite[Lemma~6.2]{PP 2017}})\\
%\end{align*}
%\end{proof}

%\vskip 2mm

%The following implies that the set of invertible elements in $\mathcal{AN}(H)$ is self-adjoint.
%
%\vskip 2mm
%
%\begin{prp}
%Let $T \in B(H)$ be invertible. Then the following is equivalent.
%
%\begin{enumerate}
%\item
%$T \in \mathcal{AN}(H)$,
%\item
%$T^* \in \mathcal{AN}(H)$.
%\end{enumerate}
%\end{prp}
%
%\begin{proof}
%
%We have only to show that the implication $(1) \rightarrow (2)$.
%

%\vskip 2mm

%\begin{rmk}
%In general, $T \in \mathcal{AN}(H)$ does not imply that $T^* \in \mathcal{AN}(H)$. (See ??\cite[Proposition~3.14]{CN 2012})
%\end{rmk}

\vskip 2mm

\begin{thm}\label{thm:functional calculas AN and AM}
Let $T\in \mathcal{B}(H)_{++}$, and $f$ be a positive continuous function 
defined on $(0,\infty)$. Then the following statements hold.
\begin{itemize}
\item[(1)] Assume that $T\in \mathcal{AN}(H)_{++}$
\begin{itemize}
\item[(i)] If $f$ is a strictly increasing function, then $f(T)\in \mathcal{AN}(H)_+$, 
\item[(ii)] if $f$ is a strictly decreasing function, then $f(T)\in \mathcal{AM}(H)_+$.
\end{itemize}
\item[(2)] Assume that $T\in \mathcal{AM}(H)_{++}$
\begin{itemize}
\item[(iii)] If $f$ is a strictly increasing function, then $f(T)\in \mathcal{AM}(H)_+$, 
\item[(iv)] if $f$ is a strictly decreasing function, then $f(T)\in \mathcal{AN}(H)_+$.
\end{itemize}
\end{itemize}
%
%In the above cases, we may assume invertibility of $T$, 
%if it is needed according to the definition of $f$.
\end{thm}

\begin{proof}
(i) has been already shown in \cite[Theorem 2.12]{GOUY 2022}. 

\noindent
Proof of (ii). Assume that $f$ is a positive and strictly decreasing function. Then $1/f(x)$ is a positive and strictly increasing function. If $T\in \mathcal{AN}(H)_{++}$. Then $f(T)^{-1}\in \mathcal{AN}(H)_{+}$ by (i). Hence $f(T)\in \mathcal{AM}(H)_{+}$ by Lemma \ref{lem:invertibility} (i).

%Put $K''=-K'$ and $F''=-F'$. Assume that $f$ is a positive continuous strictly decreasing function. By the same argument to (i), we have $K''\in \mathcal{K}(H)_+$ and $F''\in \mathcal{F}(H)_+$ such that $F''K''=0$ and $0\leq K''\leq \alpha I$. Hence $f(T)=f(\alpha)I-K''+F''\in \mathcal{AM}(H)_+$ by \cite[Theorem 5.9]{GRS2018}.

\medskip

(2) If $T\in \mathcal{AM}(H)_{++}$, then $T^{-1}\in \mathcal{AN}(H)_{++}$ by Lemma \ref{lem:invertibility} (i). 

\noindent
Proof of (iii). Assume that $f$ is a positive, continuous and strictly increasing function. Then $g(x):=f(x^{-1})$ is a positive, continuous and strictly decreasing function. Then by (ii), $f(T)=g(T^{-1})\in \mathcal{AM}(H)_+$.

(iv) can be proven by the same way to (iii). Therefore the proof is completed.
\end{proof}

Theorem \ref{thm:functional calculas AN and AM} holds for non-invertible operators, too. The proof is the same way to \cite[Theorem 2.12]{GOUY 2022}.

\begin{thm}\label{thm:arithmetric:infinity2}
For $\lambda \in (0, 1)$ and $r\in [-1,1]\setminus\{0\}$, let $f_{r,\lambda}(x)=[1-\lambda+\lambda x^r]^{1/r}$. Let $H$ be an infinite dimensional Hilbert space, and $T = U|T|\in \mathcal{B}(H)$ be  invertible and centered. If $T \in \mathcal{AN}(H)$ (resp. $T \in \mathcal{AM}(H)$). Then, there exists $\lim_{n\rightarrow\infty}\Delta_{f_{r,\lambda}}^n(T)=UL$, and  $UL\in \mathcal{AN}(H)$ (resp. $UL\in \mathcal{AM}(H)$). Moreover, the limit point does not depend on $\lambda$.
\end{thm}

\begin{proof}
By Proposition \ref{prp: power mean_centered}, we have
\begin{equation}
\Delta^{n}_{f_{r,\lambda}}(T) 
 = 
U\left[
\sum_{k=0}^{n} \left(\begin{array}{c} n \\ k \end{array}\right) a_{n,k} U^{*k}|T|^{r}U^{k}\right]^{\frac{1}{r}},
\label{eq:AN iteration}
\end{equation}
where $a_{n,k} = (1-\lambda)^k\lambda^{n-k}$.
Assume that $T\in \mathcal{AN}(H)$. Then $|T|\in \mathcal{AN}(H)_{+}$ by \cite[Corollaries 2.10 and 2.11]{NR 2019}. Moreover $|T|^{r}\in \mathcal{AN}(H)_{+}$ for $r>0$ and $|T|^{r}\in \mathcal{AM}(H)_{+}$ for $r<0$ by Theorem \ref{thm:functional calculas AN and AM}.

If $r>0$, there exist a positive number $\alpha\geq 0$, self-adjoint  $F\in \mathcal{F}(H)$ and
$K\in \mathcal{K}(H)_+$ such that 
$$ |T|^{r}=\alpha I+K+F $$
%
%$KF=0$ and $0\leq F\leq \alpha I$ 
by Theorem \ref{basic characterizations} (i).
%\cite[Theorem 2.5]{NR 2019}. 
Then by (\ref{eq:AN iteration}), we have
\begin{eqnarray*}
\Delta^{n}_{f_{r,\lambda}}(T) 
& = & 
U\left[
\sum_{k=0}^{n} \left(\begin{array}{c} n \\ k \end{array}\right) a_{n,k} U^{*k}|T|^{r}U^{k}\right]^{\frac{1}{r}}\\
& = & 
U\left[
\sum_{k=0}^{n} \left(\begin{array}{c} n \\ k \end{array}\right) a_{n,k} U^{*k}(\alpha I+K+F)U^{k}\right]^{\frac{1}{r}}\\
& = & 
U\left[\alpha I+
\sum_{k=0}^{n} \left(\begin{array}{c} n \\ k \end{array}\right) a_{n,k} U^{*k}KU^{k}
+ \sum_{k=0}^{n} \left(\begin{array}{c} n \\ k \end{array}\right) a_{n,k} U^{*k}FU^{k}\right]^{\frac{1}{r}}.
\end{eqnarray*}
From Proposition 6.1 and Corollary 6.2, there exists $K_L\in \mathcal{K}(H)_+$ and $F_L\in \mathcal{F}(H)$ such that
\begin{eqnarray*}
& & \lim_{n\to \infty}\left(
\alpha I+
\sum_{k=0}^{n} \left(\begin{array}{c} n \\ k \end{array}\right) a_{n,k} U^{*k}KU^{k}
+ \sum_{k=0}^{n} \left(\begin{array}{c} n \\ k \end{array}\right) a_{n,k} U^{*k}FU^{k}\right)\\
& = & \alpha I +K_L+F_L
=P_L\in \mathcal{AN}(H)_{+},
\end{eqnarray*}
where $P_L\in \mathcal{AN}(H)_{+}$ follows Theorem \ref{basic characterizations} (i).

Hence $L:=P_L^{\frac{1}{r}}\in \mathcal{AN}(H)_{+}$ for $r>0$
by Theorem \ref{thm:functional calculas AN and AM} (i), and there exists
$$ \lim_{n\to \infty}\Delta_{f_{r,\lambda}}^{n}(T)=UL. $$

If $r<0$, there exist a positive number $\alpha\geq 0$, $K\in \mathcal{K}(H)_+$ and $F\in \mathcal{F}(H)_+$ such that 
$$ |T|^{r}=\alpha I-K+F, $$
$KF=0$ and $0\leq K\leq \alpha I$ 
by Theorem \ref{basic characterizations} (ii).
%\cite[Theorem 5.9]{GRS2018}. 
Then by the same argument to the case $r>0$, we can prove the case $r<0$
by Theorems \ref{basic characterizations} (iii) and \ref{thm:functional calculas AN and AM} (iv). The case $T\in \mathcal{AM}(H)$ can be proven by the same way to the case $T\in \mathcal{AN}(H)$.

The limit point does not  depend on $\lambda$ by Theorem \ref{thm:convergence power mean}.
\end{proof}

\begin{rmk}
The centered condition of $T$ is not needed if $r=1$. See Remark \ref{rm:r=1 case} (i).
\end{rmk}

\vskip 2mm

\begin{qst}
Is it possible to get the same conclusion for arbitrary operator mean in Theorem~\ref{thm:arithmetric:infinity2}?
\end{qst}

\vskip 2mm

%%%%%%%%%%%%%%%%%%%%%%%%%%%%%%%%%%%%%%%%%%%%%%%%%%%%%%%%%%%%%%%%%%%%%
\section{Limit point}

In this section, we shall give
a concrete form of the
limit point of IAS 
for centered matrices.
To discuss the limit point, we recall the following result %and its proof which is 
%introduced in 
\cite[Theorem 5.1]{Y 2021}.

\begin{thm}[{\cite[Theorem 5.1]{Y 2021}}]\label{thm:iteration-finite}
Let $T\in \mathcal{M}_{m}(\C)$ 
be invertible with the polar decomposition
$T=U|T|$.
Let $\frak{a}$ be a 
non-weighted arithmetic mean.
Then the sequence $\{\Delta^n_{\frak{a}}(T)\}_{n\in \N}$ 
converges to a normal matrix $N$ 
such that
${\rm Tr}(T)={\rm Tr}(N)$ and
${\rm Tr}(|T|)={\rm Tr}(|N|)$.
\end{thm}

\begin{cor}\label{cor:limit point}
For $\lambda\in (0,1)$ and $r\in [-1,1]\setminus\{0\}$, define
$f_{r,\lambda}(x)=[1-\lambda+\lambda x^{r}]^{1/r}$. 
Let $T=U|T|\in \mathcal{M}_{m}(\C)$ 
be the polar decomposition of 
an invertible centered matrix. 
If $U=V^* diag(\alpha_{1},...,\alpha_{m})V$ 
such that $V$ is unitary and 
$\alpha_{i}\neq \alpha_{j}$ 
for $i\neq j$. Then
$$ \lim_{n\to \infty}\Delta_{f_{r,\lambda}}^{n}(T)=
U [V^*(I\circ V|T|^{r}V^*)V]^{1/r}, $$
where $X\circ Y$ means the Hadamard product of matrices $X$ and $Y$. Especially, if $r=1$, then the centered condition of $T$ can be removed.
\end{cor}

\begin{proof}
By Theorem \ref{thm:convergence power mean}, we only prove $\lambda=\frac{1}{2}$ case.
Since $\alpha_{i}\neq \alpha_{j}$ for 
$i\neq j$ and $|\alpha_{k}|=1$ for all $k=
1,2,...,m$,  
$$ \lim_{n\to \infty}
\left[\left(\frac{1+\alpha_{j}\overline{\alpha_{i}}}{2}\right)^n\right]
=I.$$
Hence by Proposition \ref{prp: power mean_centered} and the same argument as in  
\cite[Proof of Theorem 5.1]{Y 2021}, 
we have
\begin{align*}
\lim_{n\to \infty} |\Delta^n_{f_{r}}(T)|
& =
\lim_{n\to \infty}
V^{*} \left(\left[
\left(
\frac{1+\alpha_{j}\overline{\alpha_{i}}}{2}
\right)^{n}
\right]\circ V|T|^{r}V^{*}\right)^{1/r}V\\
& =
[V^{*} \left(
I\circ V|T|^{r}V^{*}\right)V]^{1/r}.
\end{align*}

If $r=1$, all of the above discussion can be done without the centered condition.
\end{proof}

\begin{exa}\label{exa:centered 3-by-3}
Let $U=\begin{pmatrix} 0 & 0 & 1 \\ 1 & 0 & 0 \\ 0 & 1 & 0 \end{pmatrix}$, 
$|T|=\begin{pmatrix} a & 0 & 0 \\ 0 & b & 0 \\ 0 & 0 & c \end{pmatrix}$ for $a,b,c>0$ and
$T=U|T|=\begin{pmatrix} 0 & 0 & c \\ a & 0 & 0 \\ 0 & b & 0 \end{pmatrix}$.
Then for each $k=1,2,...,$ 
$U^{*k}|T|U^{k}$ is a diagonal matrix, and 
$T$ is centered. For $r\in [-1,1]\setminus\{0\}$, let $f_r(x)=(\frac{1+x^r}{2})^{\frac{1}{r}}$.
In this case 
$$ N_{r}=\lim_{n\to \infty}
\Delta_{f_r}^{n}(T)=
\left(\frac{a^{r}+b^{r}+c^{r}}{3}\right)^{1/r}
\begin{pmatrix} 0 & 0 & 1 \\ 1 & 0 & 0 \\ 0 & 1 & 0 \end{pmatrix}.$$
Moreover, ${\rm Tr}(T)={\rm Tr}(N_{r})=0$,
%and 
%
\begin{align*}
\sigma(T) & = \{(abc)^{1/3},
(abc)^{1/3}\omega, (abc)^{1/3}\omega^2\}, \text{and}\\
\sigma(N_{r}) & =\{
\left(\frac{a^{r}+b^{r}+c^{r}}{3}\right)^{1/r},
\left(\frac{a^{r}+b^{r}+c^{r}}{3}\right)^{1/r}\omega, 
\left(\frac{a^{r}+b^{r}+c^{r}}{3}\right)^{1/r}\omega^2\},
\end{align*}
where $\omega\neq 1$ is a cube root of $1$.
Hence $\sigma(T)\neq \sigma(N_r)$ for $r\in [-1,1]\setminus \{0\}$.

Since $\lim_{r\to 0}(\frac{a^r+b^r+c^r}{3})^{1/r}
=(abc)^{1/3}$, we have a
limit of IAS
$$ N_{0}=\lim_{n\to \infty}\Delta^{n}(T)=
(abc)^{1/3}\begin{pmatrix} 0 & 0 & 1 \\ 1 & 0 & 0 \\ 0 & 1 & 0 \end{pmatrix}. $$
We can check that it satisfies ${\rm Tr}(T)=
{\rm Tr}(N_{0})=0$ and 
$\sigma(T)=\sigma(N_{0})$.
\end{exa}

%The matrix $T$ in Example \ref{exa:centered 3-by-3} satisfies $|T|=U^{*3}|T|U^{3}$. 
%So we 
We have the following question.

\begin{qst}\label{qst:other mean}
Can we have a concrete form of the limit of 
$\{\Delta_{\frak{m}}^{n}(T)\}$ for other 
operator mean? 
For example, consider the
logarithmic mean case 
$(f(x)=\frac{x-1}{\log x}=\int_{0}^{1}x^tdt)$.
\end{qst}

Using the matrix $T$ in Example 
\ref{exa:centered 3-by-3}, we may obtain 
a kind of the logarithmic mean for  
$3$-positive numbers as a solution of 
Question \ref{qst:other mean} if possible.
We may define any other means of
$n$-positive numbers by considering 
Question \ref{qst:other mean}.

%\section{Computation (9/24, 2022)}
%Let $a_{\lambda}(t)=1-\lambda+\lambda t$
%and $g_{\lambda}(t)=t^{\lambda}$, and let
%$T=\begin{pmatrix} 
%1 & 3 & 2 \\ 3 & 1 & 3 \\ 2 & 1 & 4 
%\end{pmatrix}$. 
%Then we compute $\Delta_{a_{0.5}}^{1000}(T)$,
%$\Delta_{a_{0.2}}^{1500}(T)$,
%$\Delta_{g_{0.5}}^{100}(T)$ and
%$\Delta_{g_{0.2}}^{100}(T)$.
%They look like limit points of iterations.
%
%\begin{itemize}
%\item[(1)] $T=\begin{pmatrix} 
%1 & 3 & 2 \\ 3 & 1 & 3 \\ 2 & 1 & 4 
%\end{pmatrix}$ and $\sigma(T)=\{6.71006, -1.89175, 1.18169\}$,
%
%\item[(2)] $\Delta_{a_{0.5}}^{1000}(T)=\begin{pmatrix} 
%0.920871 & 2.99379 & 0.526361 \\ 
%3.00265 & 1.07977 & 0.487364 \\ 
%-0.47325 & -0.539086 & 3.99936 
%\end{pmatrix}$ and \\
%$\sigma(T)=\{3.99959 + 0.716397 i, 
%3.99959 - 0.716397 i, -1.99918\}$,
%
%\item[(3)] $\Delta_{a_{0.2}}^{1500}(T)=\begin{pmatrix} 
%0.920871 & 2.99379 & 0.526361 \\ 
%3.00265 & 1.07977 & 0.487364 \\ 
%-0.47325 & -0.539086 & 3.99936 
%\end{pmatrix}$ and \\
%$\sigma(T)=\{3.99959 + 0.716397 i, 
%3.99959 - 0.716397 i, -1.99918\}$,
%
%\item[(4)] $\Delta_{g_{0.5}}^{100}(T)(T)=\begin{pmatrix} 
%0.328498 & 2.13685 & 1.68707 \\ 
%2.13685 & 0.166054 & 1.53783 \\ 
%1.68707 & 1.53783 & 5.50545 
%\end{pmatrix}$ and \\
%$\sigma(T)=\{6.71006, 
%-1.89175, 1.18169\}$,
%
%\item[(5)] $\Delta_{g_{0.2}}^{100}(T)(T)=\begin{pmatrix} 
%-0.350327 & 1.54657 & 0.278875 \\ 
%1.54657 & -0.339124 & 0.192651 \\ 
%0.278875 & 0.192651 & 6.68945 
%\end{pmatrix}$ and \\
%$\sigma(T)=\{6.71006, 
%-1.89175, 1.18169\}$.
%\end{itemize}

%%%%%%%%%%%%%%%%%%%%%%%%%%%%%%%%%%%%%%%%%%%%%%%%%%%%%%%%%%%%%%%%%%%%%%
%\section{Crossed product algebras}
\section{Iteration of induced Aluthge transformation with respect to the power mean in finite von Neumann algebras}

Let $G$ be a locally compact group, $A$ be a $C^*$-algebra and $M$ be a von Neumann algebra. 
Let $\mathrm{Aut}(A)$ denote the set of all automorphisms on $A$.

\begin{dfn}
\begin{enumerate}
\item
An action $\alpha\colon G \rightarrow \mathrm{Aut}(A)$ is said to be strongly continuous if for any $a \in A$, $\|\alpha_g(a) - a\| \rightarrow 0\ (g \rightarrow e)$.
\item
An action $\alpha\colon G \rightarrow \mathrm{Aut}(M)$ is said to be a $\sigma$-weakly continuous if for any $x \in M$ and $\phi \in M_*$, $|\phi(\alpha_g(x)) - \phi(x)| \rightarrow 0 \ (g \rightarrow e)$
\end{enumerate}

Note that if $s \in G \mapsto \alpha_s(a) \in A$ is $\sigma(A, A^*)$-continuous, then, $\alpha$ is strongly continuous.
\end{dfn}

\vskip 2mm

\begin{dfn}
Let $G$ be a locally compact group, $\mu$ be a left invariant Haar measure on $G$, and $\Delta$
be a unimodular.
A representation $\lambda\colon G \rightarrow B(L^2(G. \mu))$ is called  a left regular representation if 
for any $\xi \in L^2(G, \mu)$ and $g, h \in G$
$$
\lambda(g)\xi(h) = \xi(g^{-1}h).
$$
A representation $\rho\colon G \rightarrow B(L^2(G, \mu))$ is called a right regular representation if 
for any $\xi \in L^2(G, \mu)$ and $g, h \in G$
$$
\rho(g)\xi(h) = \Delta(t)^{\frac{1}{2}}\xi(hg).
$$
\end{dfn}

A von Neumann algebra $\mathcal{R}(G)$ generated by $\{\lambda(g)| g \in G\}$ is called a group von Neumann algebra. 
Note that $\mathcal{R}(G) = \{\lambda(g)| g \in G\}^{''}$ by von Neumann's double commutant theorem and $\rho(g) \in \mathcal{R}(G)^{'}$, the commutant of $\mathcal{R}(G)$. Note that $\mathcal{R}(G)^{'} = \{\rho(g)| g \in G\}^{''}$.

\vskip 2mm

\begin{dfn}
Let $M$ be a von Neumann algebra, G a locally compact group, and $\alpha$ $\sigma$-weakly continuous action on $M$.
Set $K = H \otimes L^2(G, \mu)$. Define a representation $\pi\colon M \rightarrow B(K)$ and a unitary representation $u\colon G \rightarrow B(K)$ by
\begin{align*}
(\pi_\alpha(x)\xi)(h) &:= \alpha_h^{-1}(x)\xi(h), \ (x \in M, \xi \in K, t \in G)\\
(u(g) \xi)(h) &:= \xi(g^{-1}h), \ (g \in G).
\end{align*}

Then, we have 
$$
\pi_\alpha((\alpha_h)(x)) = u(h)\pi_\alpha(x)u(h)^*, \ (x \in M, h \in G)
$$

The von Neumann algebra generated by $\pi_\alpha(x)$ $(x \in M)$ and $u(h)$ $(h \in G)$ is called the crossed product of $M$ and $G$, and is written as $M \rtimes_\alpha G$.
\end{dfn}

\vskip 2mm

Let $(X, \mu)$ be a probability space and $\phi$ be an invertible, measure preserving transformation of $X$. Consider the crossed product algebra $M=L^\infty(X, \mu) \rtimes_{\alpha_\phi} \Z$, where $\alpha_{\phi}(f)(x) = f(\phi(x))$ for all $f \in L^\infty(X,\mu)$ and $x \in X$. Then, there is a unitary $U$ on $ L^2(X, \mu) \otimes L^2(\Z)$ such that $\alpha_\phi(f) = UfU^* = f \circ \phi$ for all $f \in L^\infty(X, \mu)$. 
Note that linear span of  $\{U^kf| k \in \Z, f \in L^\infty(X, \mu)\}$  is strongly dense in $M$,
%Note that $\{U^kf |k \in Z, f \in L^\infty(X, \mu)\}$ is dense in $M$, 
and there is a normal state $\tau$ on $M$ such that
$$
\tau(U^kf) = \delta_{k,0}\int_Xfd\mu, \ k \in \Z, f \in L^\infty(X, \mu).
$$ 
Any element in  $\{U^kf| k \in \Z, f \in L^\infty(X, \mu)\}$  is centered because that $|U^kf| = |f| \in L^\infty(X, \mu)$, $U^{m*}|U^kf|U^m = \alpha_{-m}(|f|)$ and $U^n|U^kf|U^{n*} = \alpha_n(|f|) \in L^\infty(X, \mu)$.

\vskip 2mm

\begin{prp}
Let $T = U|T| \in M$ be centered such that $|T| \in L^\infty(X, \mu)$. 
For $r \in [-1, 1]\backslash\{0\}$,  
$\lambda\in (0,1)$ and $f_{r,\lambda}(x) =[1-\lambda+\lambda x^{r}]^{1/r}$, 
$$
\lim_{n\rightarrow\infty}\|\Delta_{f_{r, \lambda}}^n(T) - UH\|_2 = 0,
$$
where $H = (E^\phi(|T|^r))^{1/r}$.
\end{prp}

\begin{proof}
By Theorem \ref{thm:convergence power mean}, we only prove $\lambda=\frac{1}{2}$ case.

Since 
\begin{align*}
|\Delta_{f_{r, \frac{1}{2}}}^n(T)| &= \left(\frac{1}{2^n}\sum_{k=0}^n
\left(\begin{array}{c}
n\\
k
\end{array}
\right)
U^{*k}|T|^rU^k\right)^{1/r}\\
&= 
\left(\frac{1}{2^n}\sum_{k=0}^n
\left(\begin{array}{c}
n\\
k
\end{array}
\right)
|T|^r \circ \phi^k\right)^{1/r}
\end{align*}
by Proposition \ref{prp: power mean_centered},
%Theorem~\ref{thm:iteration power mean}
$|\Delta_{f_{r, \frac{1}{2}}}^n(T)|$ converges in probability  to $(E^\phi(|T|^r))^{1/r}$
by \cite[Corollary~4.4]{DS 2009}. As in the proof of \cite[Theorem~5.2]{DS 2009} we get the conclusion.
\end{proof}

\vskip 2mm

The following is the partial answer to Question 3.

\vskip 2mm

\begin{cor}\label{cor:last}
Let $T = U|T| \in M$ such that $|T| \in L^\infty(X, \mu)$. For $r \in [-1, 1]\backslash\{0\}$, $\lambda\in (0,1)$ and $f_{r,\lambda}(x) = 
[1-\lambda+\lambda x^{r}]^{1/r}$, 
$$
\Delta_{f_{r,\lambda}}^n(T) \rightarrow  UH \ (n \rightarrow\infty)
$$
in the strongly operator topology, where $H = (E^\phi(|T|^r))^{1/r}$.
\end{cor}

%%%%%%%%%%%%%%%%%%%%%%%%%%%%%%%%%%%%%%%%%%%%%%%%%%%%%%%%%
\section{Conclusion}
In this paper, we obtain  convergence properties of IASs.
To consider this matter, we  divide four cases as follows.

(i) Matrices case. In this case, for arbitrary operator mean, 
IAS of an 
invertible centered matrix converges to a normal matrix in Theorem \ref{thm:matrices case}. The limit points may depend on the 
operator mean, especially we give a concrete form of the limit point in the power mean case in Corollary \ref{cor:limit point}.
Especially, in the arithmetic mean case, the centered condition is not needed. 

(ii) Compact, $\mathcal{AN}$ and $\mathcal{AM}$ operators cases.
These operators are defined on an infinite dimensional Hilbert space.
We obtain that IAS of an 
invertible centered operator converges to a normal operator for the power mean case
in Theorems \ref{thm:convergence power mean} and  \ref{thm:arithmetric:infinity2}. 
%and \ref{thm:arithmetric:infinity}. 
%The limit point is obtained by using a conditional expectation.
It is obtained by using a kind of a generalization of von Neumann's mean ergodic theorem as in Proposition \ref{prp:extension of Dykema}.
In the arithmetic mean case, the centered condition is not needed. 

(iii) General operators case.
There is an invertible centered operator  such that IAS with respect to an arbitrary operator mean
does not converge in \cite[Theorem 5.2] {Y 2021}  (see \cite{CJL2005}). 
However, IAS with respect to an arbitrary operator mean of an invertible centered semi-hyponormal operator converges to a normal operator, and the limit point does not depend on operator means in Theorem \ref{thm:iteration general}. In the arithmetic mean case, the centered condition for an invertible semi-hyponormal operator is not needed. 

(iv) Finite von Neumann algebras case.
If $T\in M=L^\infty(X, \mu) \rtimes_{\alpha_\phi} \Z$ such that $|T|\in L^\infty(X, \mu)$, then IAS with respect to power means converges. The limit point can be represented by the conditional expectation in Corollary \ref{cor:last}.

In the all cases, all limit points do not depend on the weight of 
operator power means. In other words, we may give a classification of means by considering limit points of IASs.

%The results which we obtained can be summerized as the following table.
%
%\begin{tabular}{|c|c|c|c|c|c|c|}
% & matrix & $\mathcal{K}(H)$ & $\mathcal{AN}$ and $\mathcal{AM}$ &
% semi-hyponormal & others & finite von Neumann Algebras\\
% Arithmetic mean & converge & converge & converge & converge & not converge & converge \\ 
% Geometric mean & converge & converge & converge & converge & not converge & converge \\ 
% Power mean & converge & converge & converge & converge & not converge & converge \\ 
%arbitrary mean & converge &  &  &  & not converge &   
%\end{tabular}

\vskip 2mm
%%%%%%%%%%%%%%%%%%%%%%%%%%%%%%%%%%%%%%%%%%%%%%%%%%%%%
%\section*{Acknowledgments}
%The first author's research is supported by KAKENHI grant No. JP20K03644.

%%%%%%%%%%%%%%%%%%%%%%%%%%%%%%%%%%%%%%%%%%%%%%%%%%%%%%%%%%%%%%%

\end{document}